\documentclass[10pt,final,journal]{IEEEtran} 
\usepackage{amssymb,amsfonts}
\usepackage{algorithmic}
\usepackage{graphicx}
\usepackage{textcomp}
\def\BibTeX{{\rm B\kern-.05em{\sc i\kern-.025em b}\kern-.08em
    T\kern-.1667em\lower.7ex\hbox{E}\kern-.125emX}}
\usepackage{enumitem}
\usepackage{float}
\usepackage[tbtags]{amsmath}
\usepackage[tbtags]{mathtools}
\usepackage{accents}
\usepackage{multicol}
\usepackage[ruled,vlined]{algorithm2e}
\usepackage{cite}
\usepackage{pgfplots}

\usepackage{array,multirow}
\usepackage{hhline}
\usepackage{arydshln}

\usepackage{color,soul}

\usepackage[hidelinks=true]{hyperref}
\usepackage{cite}

\usepackage{xcolor}
\hypersetup{
    colorlinks,
    linkcolor={blue},
    citecolor={blue},
    urlcolor={blue}
}

\newtheorem{theorem}{Theorem}
\newtheorem{corollary}{Corollary}
\newtheorem{lemma}{Lemma}
\newtheorem{assumption}{Assumption}
\newtheorem{proposition}{Proposition}
\newtheorem{problem}{Problem}
\newtheorem{remark}{Remark}
\newtheorem{fact}{Fact}
\newtheorem{definition}{Definition}
\newtheorem{notation}{Notation}

\newcommand{\norm}[1]{\left\lVert#1\right\rVert}

\begin{document}
\title{On the Impulsive Formation Control of Spacecraft Under Path Constraints}
\author{Amir Shakouri
\thanks{A. Shakouri is with the Department of Aerospace Engineering, Sharif University of Technology, Tehran, Iran (e-mail: \href{mailto:a_shakouri@ae.sharif.edu}{a$\_$shakouri@ae.sharif.edu}).}}

\maketitle

\begin{abstract}
This paper deals with the impulsive formation control of spacecraft in the presence of constraints on the position vector and time. Determining a set of path constraints can increase the safety and reliability in an impulsive relative motion of spacecraft. Specially, the feasibility problem of the position norm constraints is considered in this paper. Under assumptions, it is proved that if a position vector be reachable, then the reach time and the corresponding time of impulses are unique. The trajectory boundedness of the spacecraft between adjacent impulses are analyzed using the Gerschgorin and the Rayleigh-Ritz theorems as well as a finite form of the Jensen's inequality. Some boundaries are introduced regarding the Jordan-Brouwer separation theorem which are useful in checking the satisfaction of a constraint. Two numerical examples (approximate circular formation keeping and collision-free maneuver) are solved in order to show the applications and visualize the results.
\end{abstract}

\section{Introduction}
\label{sec:I}
\IEEEPARstart{T}{he} relative spacecraft dynamics has considerably drawn the attentions due to its applications both in the formation flying and the rendezvous missions. In the simplest form, two spacecraft are considered in which a Chaser Spacecraft (CS) is the actuated system and a Target Spacecraft (TS) is located at the origin. The spacecraft are assumed to be point masses that are governed by a central gravitational force. In this context, the CS follows a path relative to the TS which is constrained dynamically and/or geometrically. Some constraints are essential to make the mission possible, while some others can be considered in order to raise the safety and reliability of the mission.

The relative motion of spacecraft can be developed as a linear time-invariant state-space model as initially proposed in \cite{in1} that is called the Clohessy-Wiltshire (CW) system. The CW model is straightforward to be implemented and optimized in the unconstrained cases \cite{in2,in3}, while the simplification assumptions are not much away from reality. However, many different models are proposed for the relative motion of spacecraft in the presence of perturbations and in the vicinity of circular or elliptical orbits \cite{in4,in4p}. 

Relative control of spacecraft is widely discussed in the literature and many different schemes are proposed. Optimal impulsive approaches based on the primer vector solutions are investigated in \cite{in5}--\cite{in7}. Gao et al. discussed the robust $H_{\infty}$ control of relative motion \cite{in8}, while solutions to the matrix inequalities are proposed by Tian and Jia \cite{in9}. Mesbahi and Hadaegh studied the formation flying control via graphs, matrix inequalities, and switching \cite{in12}. 

Constraints can be applied on the spacecraft state and/or actuator in different forms. Path constrained problems are discussed by Taur et al. \cite{in13} with impulsive time-fixed actuation. Soileau and Stern developed some necessary and sufficient conditions \cite{in14}. A method for constrained trajectory generation for micro-satellite formations is investigated by Milan et al. \cite{in15} and a finite thrust solution for state constraints is discussed by Beard and Hadaegh \cite{in16}. A model predictive control for handling the constraints is proposed by Weiss et al. \cite{in17} and Chen et al. analyzed the non-holonomic constraints \cite{in18}. For the spacecraft rendezvous with actuator saturation, a gain scheduled control is developed by Zhou et al. \cite{in19}. The global stabilization of CW system by saturated linear feedback is discussed in \cite{in20}. A covariance-based rendezvous design method is developed in \cite{in20p} by Shakouri et al. at which several constraints on the control effort, maximum impulse value, flight time, and safe zones are taken into account. 

The spacecraft actuation can be modeled by impulses in which the burning duration is negligible with respect to the time interval between adjacent burning instances. These low duration accelerations can be modeled by an impulse as a momentary change in velocity vector \cite{in5,in6,in13,in20+}. The use of hybrid impulsive-continuous actuation is studied by Sobiesiak and Damaren \cite{in21}. In \cite{in22} a multi-objective optimization is implemented by Luo et al. in order to achieve safe collision-free trajectories with admissible control efforts. 

In this paper, the impulsive behavior of the CS under equality and/or inequality path constraints are investigated. The impulse times and positions are assumed as the decision variables instead of impulse values. This point of view has advantages in considering the path constraints and disadvantages in the ignorance of the optimal solution. Under certain assumptions, it is shown that if a position vector is reachable, then the time of the next impulse and the corresponding reach time are unique (Theorem \ref{th1}). Furthermore, using the Gerschgorin circle theorem and  some results from the Rayleigh-Ritz theorem beside several spectral facts, an upper norm bound for the CS's trajectory subject to two-impulse maneuver is found (Theorem \ref{th2}) which constructs the primary contribution of the paper. Afterwards, using the Jensen's inequality, some bounding cones are introduced (Theorem \ref{th3}). The Jordan-Brouwer seperation theorem is used to show that how an area can be unreachable for the CS.
Two numerical examples are provided used to show the applications of the results. First, an approximate Circular Formation Keeping (CFK) problem is considered in which the CS finds those reachable areas as its impulse positions such that the keep-out and the approach circles constraints are not violated. It is shown how the CS can stay in an approximately circular trajectory just by two impulses. The second example is a Collision-Free Maneuver (CFM) where it is shown that how several impulse positions can be selected such that the CS do not collide with a keep-out circle with any times of impulse. This point of view can result in CFMs which are insensitive to the impulse times. 

The rest of the paper is organized as follows: First, the preliminaries regarding the notations, basic formulations, assumptions, some spectral facts, and the essential definitions are presented. Then, in Section III the main theoretical results of the paper are derived. In Section IV two numerical examples are provided in order to visualize the theoretical results. Next, some discussions are presented about the consequences and the applications. Finally, Section VII is dedicated to the concluding remarks. 

\section{Preliminaries}
\label{sec:II}
\subsection{Notation}
\label{subsec:II-A}

This subsection briefly introduces the notations that are used throughout the paper. 

Let $\mathbb{M}^{m,n}$ denote the space of $m \times n$ real (or complex) matrices and $\mathbb{M}^n$ its square analog. In addition, let $\mathbb{S}^n$ denote the space of $n$-dimensional real symmetric matrices, and $\mathbb{R}^n$ denotes the space of $n$-dimensional real vectors. The $(i,j)$th entry of a matrix $M \in \mathbb{M}^{m,n}$ is referred to by $M_{(i,j)}$ and the $i$th entry of vector $r\in \mathbb{R}^n$ is referred to by $r_{(i)}$. Upper and lower case letters are used to denote matrices and vectors, respectively. Greek letters are used to denote the scalars. For matrix $M$, we note by $M^T$ its transpose, by $M^{-1}$ its inverse (if exists), by $\textrm{null}(M)$ its nullity, and by $\textrm{rank}(M)$ its rank. 

The symbol $\|\cdot\|$ denotes a norm and specially $\|\cdot\|_p$ denotes the $p$-norm of a vector. The boundary of a set $\mathcal{S}\subset \mathbb{R}^n$ is denoted by $\partial \mathcal{S}$. We use the element-wise inequality $x\succeq 0$ to show that $x_{(i)}\geq 0$, $i=1,...,n$, and $x\succeq y$ is equivalent to $x-y\succeq0$. The symbol $\textbf{1}$ is used to denote a vector with all elements equal to 1. 

Let $\tau_1,\tau_2\in[0,\infty)$ such that $\tau_1<\tau_2$, then the following notation is used to define the time sets: 
$$\mathcal{T}_{\tau_1}=\{\tau\in \mathbb{R}|0<\tau<\tau_1\}$$
So, $\mathcal{T}_\infty=\{\tau\in \mathbb{R}|0<\tau<\infty\}$ and $\mathcal{T}_{\tau_2}-\mathcal{T}_{\tau_1}=\{\tau\in\mathbb{R}|\tau_1<\tau<\tau_2\}$. Subscript $i\in\mathbb{N}$ is used when referring to the time steps, and subscript $l\in\mathbb{N}$ is used to distinct the constraints. 

\subsection{System Model and Essentials}
\label{subsec:II-B}

In spacecraft relative motion, the Chaser Spacecraft (CS) is the actuated system and the Target Spacecraft (TS) defines the final states that needs  to be reached. Let introduce several assumptions that are made in the rest of the paper.
\begin{assumption}
\label{ass1}
Let $\alpha_{TS}\in\mathbb{R}$ denote the semi-major axis of the TS's orbit and $r\in\mathbb{R}^3$ denote the relative position of the CS with respect to the TS in an arbitrary TS-centered coordinate system. The following assumptions are made: 

\begin{enumerate}[label=\emph{(\roman*)}]
\item The two-body gravitational force is governing and no perturbations exist. 
\item The TS is in a circular orbit.
\item $\|r\|_2/\alpha_{TS}\ll1$.
\end{enumerate}
\end{assumption}

Let $r_i\in\mathbb{R}^3$ and $v_i\in \mathbb{R}^3$ denote the position and velocity vectors at step $i\in \mathbb{N}$, respectively. Suppose $r_i$ and $v_i$ are relative position and velocity of the CS defined in the RSW coordinate system of the TS (The RSW coordinate systems is defined such that its $x$-axis is in the direction of the position vector of the associated spacecraft, the $z$-axis towards the orbital angular momentum vector, and the $y$-axis completes the right-handed coordinate system). Let $t_i\in\mathcal{T}_\infty$, denote the time, $\Delta t_{i+1,i}=t_{i+1}-t_i\in\mathcal{T}_\infty$, and $F_{rr} (\cdot):\mathcal{T}_\infty\mapsto\mathbb{M}^3$ (similarly for $F_{rv}$, $F_{vr}$, and $F_{vv}$). Considering Assumption \ref{ass1} holds, the solution of the CW equations are as follows:
\begin{equation}
\label{eq:1}
r_{i+1}=F_{rr}(\Delta t_{i+1,i})r_i+F_{rv}(\Delta t_{i+1,i})v_i
\end{equation}
\begin{equation}
\label{eq:2}
v_{i+1}=F_{vr}(\Delta t_{i+1,i})r_i+F_{vv}(\Delta t_{i+1,i})v_i
\end{equation}
in which
$$F_{rr}=
\left [
  \begin{tabular}{ccc}
  $4-3 \cos{\kappa_{t}}$ & 0 & 0 \\
  $6(\sin{\kappa_{t}}-\kappa_{t})$ & 1 & 0 \\
  0 & 0 & $\cos{\kappa_{t}}$ \\
  \end{tabular}
\right ]
$$
$$
F_{rv}=\frac{1}{\kappa}
\left [
  \begin{tabular}{ccc}
  $\sin{\kappa_{t}}$ & $2(1-\cos{\kappa_{t}})$ & 0 \\
  $-2(1-\cos{\kappa_{t}})$ & $4\sin{\kappa_{t}}-3\kappa_{t}$ & 0 \\
  0 & 0 & $\sin{\kappa_{t}}$ \\
  \end{tabular}
\right ]
$$
$$
F_{vr}=\frac{d(F_{rr})}{d(\Delta t_{i+1,i})}=\kappa
\left [
  \begin{tabular}{ccc}
  $3 \sin{\kappa_{t}}$ & 0 & 0 \\
  $6(\cos{\kappa_{t}}-1)$ & 0 & 0 \\
  0 & 0 & $-\sin{\kappa_{t}}$ \\
  \end{tabular}
\right ]
$$
$$
F_{vv}=\frac{d(F_{rv})}{d(\Delta t_{i+1,i})}=
\left [
  \begin{tabular}{ccc}
  $\cos{\kappa_{t}}$ & $2\sin{\kappa_{t}}$ & 0 \\
  $-2\sin{\kappa_{t}}$ & $4\cos{\kappa_{t}}-3$ & 0 \\
  0 & 0 & $\cos{\kappa_{t}}$ \\
  \end{tabular}
\right ]
$$
where $\kappa_{t}=\kappa\Delta t_{i+1,i}$, $\kappa=\sqrt{\mu/\alpha^3_{TS}}$ is the mean motion of the TS, and $\mu$ stands for the central body (Earth) gravitational parameter. We shall use $F_{r/v,r/v}$ instead of $F_{r/v,r/v}(\cdot)$ for simplicity. Here, another assumption is made to avoid singularities in the matrices defined above. 

\begin{assumption}
\label{ass2}
The time interval between any two subsequent impulses should be less than $\pi/\kappa$, i.e., $\forall i\in\{1,...,n-1\}$, $\Delta t_{i+1,i}\in\mathcal{T}_{\pi/\kappa}$. 
\end{assumption}

It is worth mentioning that for a relative motion of spacecraft with a total flight time of $\Delta t_{total}\in[0,\infty)$, if the number of impulses satisfy $n\geq \min\{\text{floor}[\Delta t_{total}/(\pi/\kappa)]+1,2\}$, then Assumption \ref{ass2} can be satisfied.

Let us introduce the position and velocity vectors, $r_i^-$ and $v_i^-$, in which define the position and velocity vectors before applying an impulse vector, $\Delta v_i\in\mathbb{R}^3$, at $i$th step. After applying the impulse vector, the position and velocity vectors become $r_i^+=r_i^-\equiv r_i$ and $v_i^+=v_i^-+\Delta v_i$. So, the position and velocity after $\Delta t_{i+1,i}$ are found from \eqref{eq:1} as
\begin{equation}
\label{eq:3}
r_{i+1}=F_{rr}r_i+F_{rv}v_i^+=F_{rr}r_i+F_{rv}(v_i^-+\Delta v_i)
\end{equation}
\begin{equation}
\label{eq:4}
v_{i+1}=F_{vr}r_i+F_{vv}v_i^+=F_{vr}r_i+F_{vv}(v_i^-+\Delta v_i)
\end{equation}
Therefore, form \eqref{eq:2} knowing $r_i$, $r_{i+1}$, and $v_i^-$, the impulse vector $\Delta v_i$ is
\begin{equation}
\label{eq:5}
\Delta v_i=F^{-1}_{rv}(r_{i+1}-F_{rr}r_i)-v_i^-
\end{equation}

\begin{remark}
\label{rem1}
For an $n$-impulse relative motion between known $r_1$, $r_n$, $v_1^-$, and $v_n^+$, which is usually the case, the decision variables can be chosen to be one of the following sets:
\begin{enumerate}[label=\emph{(\roman*)}]
\item $\{r_i\in\mathbb{R}^3|i=2,...,n-1\}\cup\{\Delta t_{i+1,i}\in\mathcal{T}_{\pi/\kappa}|i=1,...,n-1\}$
\item $\{\Delta v_i\in\mathbb{R}^3|i=1,...,n-2\}\cup\{\Delta t_{i+1,i}\in\mathcal{T}_{\pi/\kappa}|i=1,...,n-1\}$
\item $\{r_i\in\mathbb{R}^3|i=2,...,n-1\}\cup\{\|\Delta v_i\|_2\in\mathbb{R}^3|i=1,...,n-1\}$
\end{enumerate}
\end{remark}
Each of the above sets has $2n-3$ members. So, for an $n$-impulse trajectory, a number of $2n-3$ decision variables ($n-2$ vectors and $n-1$ scalars) are needed to determine the whole trajectory. 

Consider the following matrices:
\begin{equation}
\label{eq:7}
F_2(t,\tau)=F_{rv}(t)F_{rv}^{-1}(\tau)
\end{equation}
\begin{equation}
\label{eq:8}
F_1(t,\tau)=F_{rr}(t)-F_2(t,\tau)F_{rr}(\tau)
\end{equation}
in which $\tau\in\mathcal{T}_{\pi/\kappa}$ and $t\in\mathcal{T}_{\tau}$. Suppose $\tau$ to be fixed, so we use $F_1(t,\tau)\equiv F_1(t)$ and $F_2(t,\tau)\equiv F_2(t)$ for simplicity. From \eqref{eq:1} to \eqref{eq:4} it can be shown that how the CS's trajectory behaves in the time domain subjected to fixed initial ($r_i$) and final positions ($r_{i+1}$) as well as the total flight time ($\Delta t_{i+1,i}$). Using the forms defined in \eqref{eq:7} and \eqref{eq:8}:
\begin{equation}
\label{eq:9}
r(t,\Delta t_{i+1,i})=F_1(t,\Delta t_{i+1,i})r_i+F_2(t,\Delta t_{i+1,i})r_{i+1}
\end{equation}
in which we use simply $r\equiv r(t,\Delta t_{i+1,i})$ for fixed $\Delta t_{i+1,i}$. 

\subsection{Spectral Analysis}
\label{subsec:II-C}

The spectral properties of some matrices are needed to be analyzed to be used further in obtaining the results. Let $\lambda_{11}(t),\lambda_{22}(t)\in\mathbb{R}^3$ denote the eigenvalues of  $F_1^T(t)F_1(t)$ and $F_2^T(t)F_2(t)$, respectively, that are sorted in vectors such that ${\lambda_{ii}}_{(j)}(t)\leq{\lambda_{ii}}_{(j+1)}(t)$ for $i,j\in\{1,2\}$. Now, consider the following $6\times6$ symmetric block form matrix: 
\begin{equation}
\label{eq:10}
\widehat{F}(t)=
\left [
  \begin{tabular}{cc}
  $F_1^T(t)F_1(t)$ & $F_1^T(t)F_2(t)$ \\
  $F_2^T(t)F_1(t)$ & $F_2^T(t)F_2(t)$ \\
  \end{tabular}
\right ]
\end{equation}
Denote the eigenvalues of $\widehat{F}(t)$ by $\hat{\lambda}_{(1)}(t)\leq\hat{\lambda}_{(2)}(t)\leq...\leq\hat{\lambda}_{(6)}(t)$. The following facts are analytically or numerically evaluated: 

\begin{fact}
\label{f1}
Under Assumptions \ref{ass1} and \ref{ass2} the following statements hold for matrices $F_1$ and $F_2$:
\begin{enumerate}[label=\emph{(\roman*)}]
\item $\text{rank}[F_1^T(t)F_1(t)]=3$, $\text{rank}[F_2^T(t)F_2(t)]=3$, and both have three real nonzero eigenvalues.
\item At any $t\in\mathcal{T}_{\tau}$, the following property holds for the eigenvalues of $F_1^TF_1$ and $F_2^TF_2$: 
\begin{equation}
\label{eq:11}
{\lambda_{11}}(t)-{\lambda_{22}}(\tau-t)=0
\end{equation}
\item For $i\in\{1,2\}$, at any $t\in\mathcal{T}_{\tau}$, ${\lambda_{ii}}_{(3)}<1$ for $0<\tau<(\pi/\kappa)/2$, ${\lambda_{ii}}_{(3)}=1$ for $\tau=(\pi/\kappa)/2$, and ${\lambda_{ii}}_{(3)}>1$ for $(\pi/\kappa)/2<\tau<\pi/\kappa$. 
\end{enumerate}
\end{fact}

\begin{fact}
\label{f2}
Under Assumptions \ref{ass1} and \ref{ass2} the following statements hold for matrix $\widehat{F}$:
\begin{enumerate}[label=\emph{(\roman*)}] 
\item $\text{rank}[\widehat{F}(t)]=3$, $\text{null}[\widehat{F}(t)]=3$, and it has three eigenvalues of zero, i.e., $\hat{\lambda}_{(1)}(t)=\hat{\lambda}_{(2)}(t)=\hat{\lambda}_{(3)}(t)=0$.
\item At any $t\in\mathcal{T}_{\tau}$, the entries of $\widehat{F}(t)$ has the following property for $i=1,2,3$: 
\begin{equation}
\label{eq:12}
\sum_{j=1}^6\|\widehat{F}_{(i,j)}(t)\|_1-\|\widehat{F}_{(i+3,j)}(\tau-t)\|_1=0
\end{equation}
\item The eigenvalues of $\widehat{F}(t)$ are constant over time or they have a single extremum at $t=\tau/2$. 
\item At any $t\in\mathcal{T}_{\tau}$, $\hat{\lambda}_{(6)}<1$ for $0<\tau<(\pi/\kappa)/2$, $\hat{\lambda}_{(6)}=1$ for $\tau=(\pi/\kappa)/2$, and $\hat{\lambda}_{(6)}>1$ for $(\pi/\kappa)/2<\tau<\pi/\kappa$. 
\end{enumerate}
\end{fact}

\subsection{Definitions}
\label{subsec:II-D}

\begin{notation}
Let $\Gamma_i^j(n)$ denote a trajectory such that $n$ impulses are used starting from index $i$ and ending in $j$ such that the first and the last impulses are applied at $i$ and $j$, respectively. For example, $\Gamma_1^3(2)$ denotes a two-impulse trajectory that starts from $r_1$, $v_1^-$ at $t_1$, and ends in $r_3$ at $t_3$. 
\end{notation}

It should be noted that the symbol $\Gamma_i^j(n)$ do not give any knowledge about the decision variables and the initial/final states of the trajectory. So, $\Gamma_i^j(n)$ alone cannot define a relative spacecraft motion trajectory even for two-impulse missions. 

\begin{definition}
\label{def1}
Let $r\in\mathbb{R}^3$ denote the position vector at $t\in\mathcal{T}_\infty$. Consider the sets $\widetilde{\mathcal{T}}_l\subset\mathcal{T}_\infty$ and $\widetilde{\mathcal{R}}_l\subset\mathbb{R}^3$ at $l\in\{1,2,...,m\}$ for $m\in\mathbb{N}$. Path Constraints (PCs) are those constraints that can be stated as follows: 
\begin{equation}
\label{eq:13}
\forall t\in\widetilde{\mathcal{T}}_l\text{ : }r\in\widetilde{\mathcal{R}}_l
\end{equation}
\end{definition}

In this paper we are dealing with a special kind of PCs. Suppose $\tilde{r}_l\in\mathbb{R}^3$, $\rho_l^{\prime},\rho_l^{\prime\prime}\in\mathbb{R}$, and $\tilde{t}_l\in\mathcal{T}_{\infty}$ are predefined parameters at $l=1,...,m$. The general PCs defined in \eqref{eq:13} can be reduced to an inequality form such that:
\begin{equation}
\label{eq:14}
\widetilde{\mathcal{T}}_l=\mathcal{T}_{\tilde{t}_l},
\quad
\widetilde{\mathcal{R}}_l=\left\{r\in\mathbb{R}^3|\rho_l^{\prime}\leq\|r-\tilde{r}_l\|_2\leq \rho_l^{\prime\prime}\right\}
\end{equation}

Each PC of the form \eqref{eq:14} restricts the CS's path between two spheres at a time interval. The form of inequality PCs introduced in \eqref{eq:14} can be used for trajectories that avoid collisions independent from the transfer time. The following feasibility problem, states the main subject of the paper.

\begin{problem}
Find the decision variables (discussed in Remark \ref{rem1}) such that satisfies a PC of the form \eqref{eq:14}. 
\end{problem}

\begin{remark}
\label{rem2}
In \eqref{eq:14}, assuming $\rho_l^{\prime\prime}\rightarrow\infty$, the PC defines a forbidden region at $t\in\mathcal{T}_{\tilde{t}_l}$. This region is bounded by $\|r-\tilde{r}_l\|_2=\rho_l^{\prime}$. This kind of PCs can be used to define collision-free trajectories that are robust with respect to actuator fault and failure. 
\end{remark}

The inequality form of \eqref{eq:14} can turns to equality if $\widetilde{\mathcal{T}}_l=\{\tilde{t}_l\}$ and $\rho_l^{\prime}=\rho_l^{\prime\prime}=0$. So, an equality PC can be stated in the following form: 

\begin{equation}
\label{eq:15}
\widetilde{\mathcal{T}}_l=\{\tilde{t}_l\},
\quad
\widetilde{\mathcal{R}}_l=\{\tilde{r}_l\}
\end{equation}

\begin{remark}
\label{rem3}
A spacecraft trajectory with known initial and final positions, $r_1$ and $r_n$, which is usually the case, essentially has two PCs of the equality form; the first is $\|r-\tilde{r}_1\|_2=0$ at $t=0$ and the second is $\|r-\tilde{r}_n\|_2=0$ at $t=t_n$. 
\end{remark}

A two-point constrained single-impulse reachability (or simply ``reachability'') can be defined in the context of this paper which is an impulsive reachability that is constrained in order to achieve an initial and a final position. 

\begin{definition}
\label{def2}
Let $t_j-t_i\in\mathcal{T}_{\pi/\kappa}$ and $\tau\in [t_i,t_j]$. A position vector $\tilde{r}$ is called reachable in $\Gamma_i^j(2)$ at $t\in[t_i,\tau]\subseteq[t_i,t_j]$, if there exists a $Δv_i\in\mathbb{R}^3$ such that $r=\tilde{r}$ at $t\in[t_i,\tau]$ subject to $r=r_i$ at $t=t_i$ and $r=r_j$ at $t=t_j$. A position vector that is not reachable is called unreachable. 
\end{definition}

\subsection{Time Uniqueness}
\label{subsec:III-A}

In this subsection it is proved that if a point in the space under Assumptions \ref{ass1} and \ref{ass2} is reachable, so the corresponding total time of flight and the current time are unique. Consider the following results: 

\begin{lemma}
\label{lem1}
Let $\mathfrak{T}^{\prime},\mathfrak{T}^{\prime\prime}\subset\mathcal{T}_{\pi/\kappa}$ and $\mathfrak{R}_r,\mathfrak{R}_v\subset \mathbb{R}^3$, then functions $v_t(\cdot):\mathfrak{T}^{\prime}\mapsto\mathfrak{R}_v$, $r_v(\cdot):\mathfrak{T}^{\prime\prime}\times\mathfrak{R}_v\mapsto\mathfrak{R}_r$, and $r_t(\cdot):\mathfrak{T}^{\prime}\times\mathfrak{T}^{\prime\prime}\mapsto\mathfrak{R}_r$ are injective such that: 
\begin{equation}
\label{eq:16}
v_t(t^{\prime})=F_{rv}^{-1}(t^{\prime})\left[r_j-F_{rr}(t^{\prime})r_i\right]-v_i^-
\end{equation}
\begin{equation}
\label{eq:17}
r_v(\Delta v_i,t^{\prime\prime})=F_{rr}(t^{\prime\prime})r_i+F_{rv}(t^{\prime\prime})v_i^-+F_{rv}(t^{\prime\prime})\Delta v_i
\end{equation}
\begin{equation}
\label{eq:18}
\begin{split}
r_t(t^{\prime},t^{\prime\prime})=\left[F_{rr}(t^{\prime\prime})-F_{rv}(t^{\prime\prime})F_{rv}^{-1}(t^{\prime})F_{rr}(t^{\prime})\right]r_i \\
+F_{rv}(t^{\prime\prime})F^{-1}_{rv}(t^{\prime})r_j
\end{split}
\end{equation}
\end{lemma}
\begin{IEEEproof}
First it should be noted that $F_{rr}(t^{\prime})$ and $F_{rv}(t^{\prime})$ are injective maps from $\mathfrak{T}$ to some subsets of $\mathbb{M}^3$. Function $v_t(t^{\prime})$ generates an impulse vector to transfer the CS from $r_i$ and $v_i^-$ at $t_i$ to $r_j$ at $t_j$. From the physics of the problem, obviously, each impulse ($\forall t^{\prime}\in\mathcal{T}_\infty$) has a unique time to transfer and \eqref{eq:16} is injective. To show the uniqueness of \eqref{eq:17} consider $r_v(\Delta\underaccent{\bar}{v}_i)=r_v(\Delta\bar{v}_i)$ that leads to $F_{rv}(t^{\prime\prime})(\Delta\underaccent{\bar}{v}_i-\Delta\bar{v}_i)=0$. Given that $\forall t^{\prime\prime}\in\mathcal{T}_{\pi/\kappa}$, $\textrm{rank}⁡[F_{rv}(t^{\prime\prime})]=3$, according to the rank-nullity theorem $\textrm{null}⁡[F_{rv}(t^{\prime\prime})]=3$. So, the only solution to $F_{rv}(t^{\prime\prime})(\Delta\underaccent{\bar}{v}_i-\Delta\bar{v}_i)=0$ is $\Delta\underaccent{\bar}{v}_i=\Delta\bar{v}_i$. Function $r_t(\cdot)$ is $r_v(\cdot)$ composed with the injective function $\accentset{\star}{v}_t(t^{\prime},t^{\prime\prime})=[v_t^T(t^{\prime})\quad t^{\prime\prime}]^T$, i.e., $r_t=r_v\circ \accentset{\star}{v}_t$. The composition of injective functions is injective, so $r_t(\cdot)$ is injective.
\end{IEEEproof}
\begin{lemma}
\label{lem2}
Let $i,j\in\mathbb{N}$, $\mathcal{R}_{i,j}(t)\subset\mathbb{R}^3$, and $t_i,t_j,t\in\mathcal{T}_{\infty}$ where $t_i<t<t_j$ and $\Delta t_{j,i}\in\mathcal{T}_{\pi/\kappa}+\mathcal{T}_{t_i}-\mathcal{T}_{t}$, such that
\begin{equation}
\label{eq:19}
\begin{split}
\mathcal{R}_{i,j}(t)=\left\{r\in\mathbb{R}^3|\forall \Delta t_{j,i}\in\mathcal{T}_{\pi/\kappa}+\mathcal{T}_{t_i}-\mathcal{T}_{t}\text{ : }r=\right. \\
\left.F_1(t-t_i,\Delta t_{j,i})r_i+F_2(t-t_i,\Delta t_{j,i})r_j\right\}
\end{split}
\end{equation}
then $r_j$ is reachable in $\Gamma_i^j(2)$ at $t\in[t_i,t_j]$ if and only if $r\in\mathcal{R}_{i,j}(t)$
\end{lemma}
\begin{IEEEproof}
It can be directly concluded from \eqref{eq:9} that if $r\in\mathcal{R}_{i,j}(t)$ then $\exists\Delta t_{i,j}\in\mathcal{T}_{\pi/\kappa}$ such that the CS reaches $r$ at $t$. For the other direction, we know that if $r\notin\mathcal{R}_{i,j}(t)$ then the CS cannot reach $r$ at $t$, $\forall\Delta t_{i,j}\in\mathcal{T}_{\pi/\kappa}$. 
\end{IEEEproof}
\begin{remark}
\label{rem5}
In Lemma \ref{lem2}, the elements of set $\mathcal{R}_{i,j}(t)$ are the outputs of a nonlinear mapping from $\mathcal{T}_{\pi/\kappa}+\mathcal{T}_{t_i}-\mathcal{T}_{t}\subset\mathbb{R}$ to a higher dimension in $\mathbb{R}^3$. In fact, the locus of the elements of $\mathcal{R}_{i,j}(t)$ is a three-dimensional curve that starts from $r\to r_j$ (at $\Delta t_{j,i}\to t$) and ends in some infinite $r$ (at $\Delta t_{j,i}\to \pi/\kappa$). 
\end{remark}

\begin{proposition}
\label{prop1}
Let $\tau_1,\tau_2\in\mathcal{T}_{\pi/\kappa}$, then $\mathcal{R}_{1,2}(\tau_1 )\cap\mathcal{R}_{1,2}(\tau_2 )\ne\varnothing$ if and only if $\tau_1=\tau_2$.
\end{proposition}
\begin{IEEEproof}
The necessity is obvious. For sufficiency, a proof by contradiction is used. Suppose $\tau_1\ne\tau_2$ and $\mathcal{R}_{1,2}(\tau_1)\cap\mathcal{R}_{1,2}(\tau_2)\ne\varnothing$, then there exists a position $r_t$ and two time intervals $\Delta t_{21}^\prime$ and $\Delta t_{21}^{\prime\prime}$ such that $r_t(\Delta t_{21}^\prime,\tau)=r_t(\Delta t_{21}^{\prime\prime},\tau)$. Thus, from Lemma \ref{lem1} if $\Delta t_{21}^{\prime}\ne\Delta t_{21}^{\prime\prime}$, then no equality exists and if $\Delta t_{21}^{\prime}=\Delta t_{21}^{\prime\prime}$, then $\tau_1=\tau_2$ which demonstrates a contradiction.
\end{IEEEproof}
\begin{theorem}
\label{th1}
Suppose Assumptions \ref{ass1} and \ref{ass2} hold. If the position vector $r$ is reachable in $\Gamma_i^{i+1}(2)$ at $t\in[t_i,t_{i+1}]$, then the corresponding $t$ and $\Delta t_{i+1,i}$ are unique.
\end{theorem}
\begin{IEEEproof}
From the result of Lemma \ref{lem2}, the reachability of $r$ in $\Gamma_i^{i+1}(2)$ at $t\in[t_i,t_{i+1}]$ is guarantied if and only if there exist a $t$ and $\Delta t_{i+1,i}$ such that $r=r_t(\Delta t_{i+1,i},t)$. From Lemma \ref{lem1}, $r_t(\cdot)$ is injective. Therefore, having the output, the corresponding inputs are unique.
\end{IEEEproof}

\textcolor{black}{Theorem \ref{th1} individually concludes that a two-impulse rendezvous maneuver which is restricted to reach a point in the space (except of initial and final locations), has a unique solution. This can be used for cases that we need to hit a target between initial and final locations. Furthermore, Theorem \ref{th1} can be used for an initial relative orbit determination using three vectors (similar to the Gibbs method in two-body problem) which can be a subject for the future researches.}

\textcolor{black}{Regarding Theorem \ref{th1}, the reader can refer to the work done by Wen et al. in \cite{in23} at which the reachability problem of impulsive maneuvers for nonlinear unperturbed problems is investigated through analytic geometry.}

\section{Main Results}
\label{sec:III}

\subsection{Trajectory Boundedness}
\label{subsec:III-B}

In this subsection an upper bound on the CS's trajectory subjected to fixed initial and final positions are presented. It is shown that for a $2$-impulse transfer between $r_i$ at $t_i$ and $r_{i+1}$ at $t_{i+1}$ the trajectory has upper bounds on the norm of the position vector, depending on the initial ($r_i$) and final positions ($r_{i+1}$) as well as the total flight time ($\Delta t_{i+1,i}$). The results of this section can be directly used for the design of the decision variables in a constrained formation of spacecraft. First we need two lemmas: 

\begin{lemma}[Gerschgorin]
\label{lem3}
Let $M\in\mathbb{M}^n$, with associated eigenvalues of $\mu_i$, $i=1,...,n$, and let
\begin{equation}
\label{eq:20}
\rho_i(M)=\sum_{j=1,j\neq i}^n\|M_{(i,j)}\|_1,\quad 1\leq i\leq n
\end{equation}
denote the deleted absolute row sums of $M$. Then all the eigenvalues of $M$ are located in the union of $n$ circles (i.e., $\forall i\in{1,...,n}:\mu_i\in \mathcal{G}(M)$): 
\begin{equation}
\label{eq:21}
\mathcal{G}(M)=\bigcup_{i=1}^n\{z\in\mathbb{R}|\|z-M_{(i,i)}\|_1\leq \rho_i(M)\}
\end{equation}
\end{lemma}
\begin{IEEEproof}
A proof can be found in \cite{b1}, pp. 344. 
\end{IEEEproof}
\begin{lemma}[Rayleigh-Ritz]
\label{lem4}
Let $M\in\mathbb{M}^n$ be Hermitian and $x\in\mathbb{R}^n$. Let $\mu_{max}$ and $\mu_{min}$ denote the maximum and the minimum eigenvalue of $M$, respectively. Then
\begin{equation}
\label{eq:22}
\mu_{min}\|x\|_2^2\leq x^TMx\leq\mu_{max}\|x\|_2^2
\end{equation}
\end{lemma}
\begin{IEEEproof}
A proof can be found in \cite{b1}, pp. 176. 
\end{IEEEproof}
\begin{theorem}
\label{th2}
Suppose Assumptions \ref{ass1} and \ref{ass2} hold. The $2$-impulse  trajectory of the CS subject to $r_i$ at $t_i$ and $r_{i+1}$ at $t_{i+1}$ is bounded by a sphere centered at the origin with a radius of $\delta_{i+1,i}$, i.e., $\forall t\in\mathcal{T}_{t_{i+1}}-\mathcal{T}_{t_{i}}:\|r\|_2\leq\delta_{i+1,i}$, such that:
\begin{equation}
\label{eq:23}
\delta_{i+1,i}=\sigma(\Delta t_{i+1,i})\sqrt{\|r_i\|_2^2+\|r_{i+1}\|_2^2}
\end{equation}
\begin{equation}
\label{eq:24}
\sigma(t)=
\left\{
\begin{tabular}{cc}
$1 \hfill$ & $t\in(0,0.5\pi/\kappa]\hfill$ \\
$0.5\sqrt{2}\sec(0.5\kappa t)$ & $t\in[0.5\pi/\kappa,\pi/\kappa)$
\end{tabular}
\right.
\end{equation}
\end{theorem}
\begin{IEEEproof}
The two-impulse trajectory subjected to fixed initial ($r_i$) and final ($r_{i+1}$) positions with a fixed time of flight ($\Delta t_{i+1,i}$), has a position vector $r$ at $t$ that is introduced previously in \eqref{eq:11}. From \eqref{eq:11}, the squared 2-norm of $r$ can be written as
\begin{equation}
\label{eq:25}
\|r\|_2^2=r^Tr=\begin{bmatrix} r_i \\ r_{i+1} \end{bmatrix}^T
\widehat{F}
\begin{bmatrix} r_i \\ r_{i+1} \end{bmatrix}
\end{equation}

From Lemma \ref{lem4} the above identity has the following upper bound: 
\begin{equation}
\label{eq:26}
\lambda_{max}(t)\norm{\begin{bmatrix} r_i \\ r_{i+1} \end{bmatrix}}_2^2=\lambda_{max}(t)(\|r_i\|_2^2+\|r_{i+1}\|_2^2)
\end{equation}
in which $\lambda_{max}(t)$ is the maximum eigenvalue of $\widehat{F}$ at $t\in\mathcal{T}_{\Delta t_{i+1,i}}$. Let us determine an upper bound for $\lambda_{max}(t)$ in order to make the formulas independent from time. Suppose:
\begin{equation}
\label{eq:27}
\sigma^2=\sup_{t\in\mathcal{T}_{\Delta t_{i+1,i}}}\lambda_{max}(t)
\end{equation}
So, using \eqref{eq:25} to \eqref{eq:27}:
\begin{equation}
\label{eq:28}
\|r\|_2^2\leq\sigma^2(\|r_i\|_2^2+\|r_{i+1}\|_2^2)=\delta_{i+1,i}^2
\end{equation}

From Facts \ref{f1} and \ref{f2} it can be concluded that $\sigma^2=1$ at $\Delta t_{i+1,i}\leq(\pi/\kappa)/2$, since at the initial and final positions all of the positive eigenvalues reach unity which is the solution of \eqref{eq:27}. At $(\pi/\kappa)/2<\Delta t_{i+1,i}<\pi/\kappa$, $\sigma^2>1$ and is equal to its extremum value at $t=\Delta t_{i+1,i}/2$. Lemma \ref{lem3} provides an upper bound on the eigenvalues at $t=\Delta t_{i+1,i}/2$ in which from Fact \ref{f2}, the upper boundary of Gershgorin circle has also an extremum at $t=\Delta t_{i+1,i}/2$ and is equal to the following equation: 
\begin{equation}
\label{eq:29}
\begin{split}
\max_{k=1,...,6}\left\{\widehat{F}_{(k,k)}(\Delta t_{i+1,i}/2)+\rho_k\left[\widehat{F}(\Delta t_{i+1,i}/2)\right]\right\} \\
=\frac{1}{2}\sec^2\left(\frac{\kappa \Delta t_{i+1,i}}{2}\right)
\end{split}
\end{equation}
Thus, at $(\pi/\kappa)/2<\Delta t_{i+1,i}<\pi/\kappa$, $\sigma(\Delta t_{i+1,i})=0.5\sqrt{2}\sec(0.5\kappa \Delta t_{i+1,i})$, and the theorem can be proved considering this result by taking the square root of \eqref{eq:28}. 
\end{IEEEproof}

\textcolor{black}{Eq. \eqref{eq:24} of Theorem \ref{th2} checks that even if the true anomaly of the TS has been changed less than $90^\circ$ or more. Note that the TS can be just a reference for the coordination of the CS and can be assumed a virtual point of reference. Thus, without loss of generality, one can assume $\|r_{i+1}\|=0$ in which concludes that for a spacecraft rendezvous under Assumption \ref{ass1}, if the decision variable $\Delta t_{i+1,i}$ be considered less than $0.5\pi/\kappa$, then the spacecraft will not increase its distance from the destination. If Assumption \ref{ass1} approximately holds for real case applications, this kind of maneuver can assure that the approximation error would not grow.}

\begin{remark}
\label{rem6}
The trajectory bound of Theorem \ref{th2} becomes more tight by decreasing the norm of the initial and final positions (i.e., making the impulses closer to the TS at the origin). The bound is time-independent and for $\Delta t_{i+1,i}\leq(\pi/\kappa)/2$ it is also independent from the total flight time. For $(\pi/\kappa)/2<\Delta t_{i+1,i}<\pi/\kappa$ the bound increases by increasing the total flight time and approaches infinity at $\Delta t_{i+1,i}\rightarrow\pi/\kappa$.
\end{remark}

\begin{figure}[!h]
\centering\includegraphics[width=0.5\linewidth]{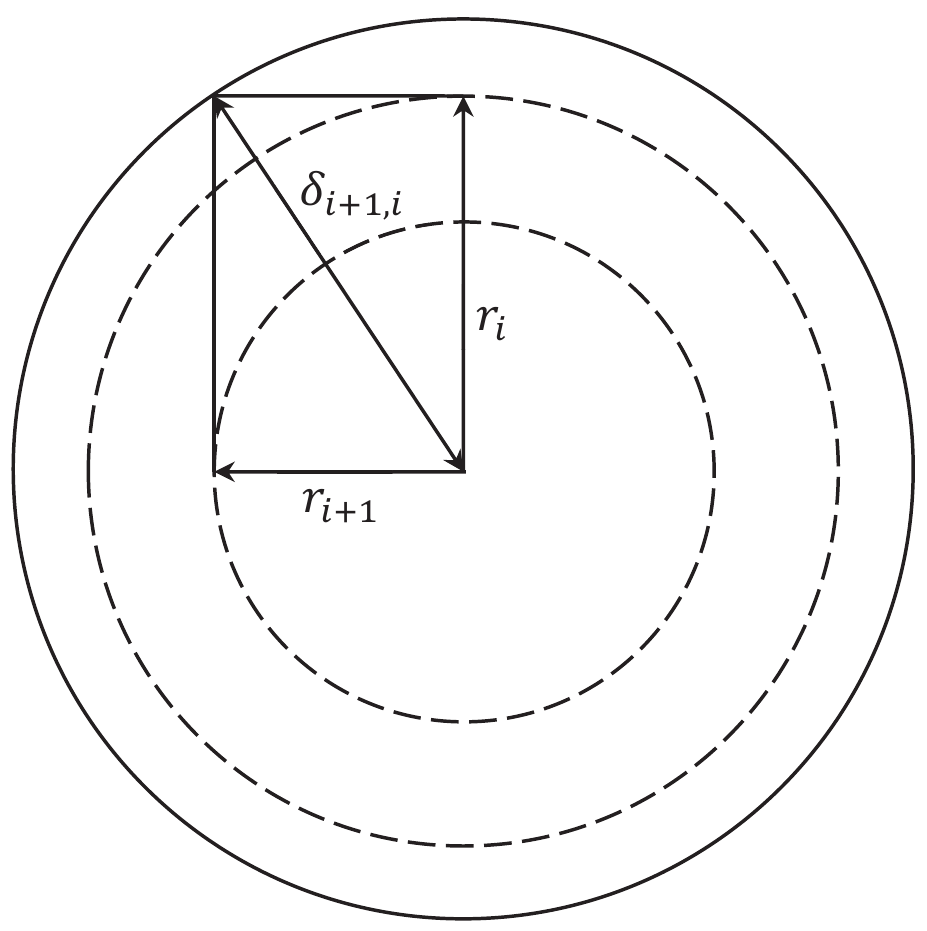}
\caption{The dashed spheres show the location of CS at steps $i$ and $i+1$. The solid sphere shows the boundary of CS's trajectory for $\Delta t_{i+1,i}\leq(\pi/\kappa)/2$.}
\label{fig:2}
\end{figure}

\begin{figure}[!h]
\centering\includegraphics[width=0.7\linewidth]{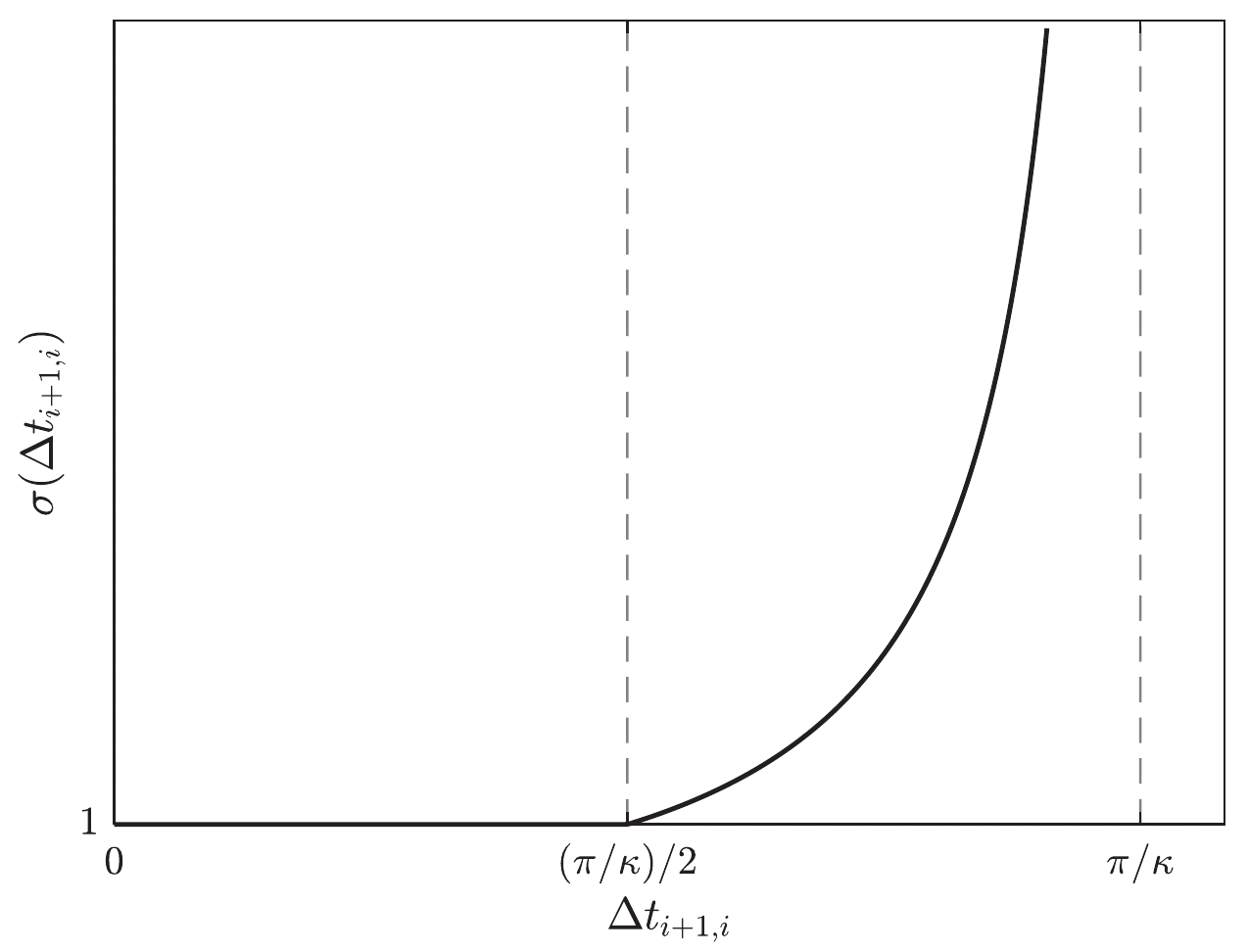}
\caption{Behavior of $\sigma$ as a function of $\Delta t_{i+1,i}$.}
\label{fig:3}
\end{figure}

A geometric interpretation of the 2-norm bounds defined in Theorem \ref{th2} is illustrated is Fig. \ref{fig:2}. This figure shows the bound for $\Delta t_{i+1,i}\leq(\pi/\kappa)/2$. For $\Delta t_{i+1,i}>(\pi/\kappa)/2$ the solid sphere becomes bigger. Fig. \ref{fig:3} shows the behavior of $\sigma$ as a function of $\Delta t_{i+1,i}$. 


Theorem \ref{th2} provides an upper bound for the 2-norm of the CS's trajectory  which is useful in the design of the impulse positions.  It can be concluded that for a multiple impulse relative motion with uncertain impulse positions, if it is known that each impulse position at $i$ is bounded by a sphere with a radius of $\rho_i$, and the time interval between adjacent impulses be bounded as $\Delta t_{i+1,i}<(\pi/\kappa)/2$, then the whole trajectory is bounded in a sphere with a radius of
\begin{equation}
\label{eq:30}
\rho=\sqrt{2}(\max_{i}\{\rho_i\})
\end{equation}

The CS's trajectory can be bounded by the use of a cone, i.e., some bounds for the scalar parameter $s^Tr$ which $s\in\mathbb{R}^3$ is a known and probably fixed (and unit) vector. Both upper and lower bounds for $s^Tr$ can restrict the CS's trajectory to lie inside or outside a cone with the apex on the origin. First, consider a special finite form of the Jensen's inequality \cite{b2} in which the proof is omitted: 

\begin{lemma}[Jensen]
\label{lem5}
Let $s,r\in\mathbb{R}^n$ and $s\succeq0$. Consider a real function $f:\mathbb{R}\mapsto\mathbb{R}$. If $f$ is convex \cite{b3} then 
\begin{equation}
\label{eq:31}
f\left( \frac{s^Tr}{\textbf{1}^Ts}\right) \leq\frac{\sum s_{(i)} f(r_{(i)})}{\textbf{1}^Ts}
\end{equation}
and if $f$ is concave ($-f$ is convex) then the above inequality holds with a change in its direction.
\end{lemma}

From the Jensen's inequality in Lemma \ref{lem5}, upper and/or lower bounds can be found for $s^Tr$ by choosing any convex (or concave) function $f$. Suppose $f$ be a convex nondecreasing invertible function (for example $f(x)=y^x$, $y\geq1$). Then inequality \eqref{eq:31} simplifies to
\begin{equation}
\label{eq:32}
s^Tr\leq(\textbf{1}^Ts)f^{-1}\left[ \frac{\sum s_{(i)} f(r_{(i)})}{\textbf{1}^Ts}\right] 
\end{equation}
Inequality \eqref{eq:32} defines an upper bound for $s^Tr$, if some priory knowledge about the $r$ exists. For example, suppose it is known that the elements of $r$ are bounded above, i.e., $r\preceq r^{+}$. So, from the nondecreasing property of $f$ we have $f(r_{(i)})\leq f(r^{+}_{(i)})$. Thus,
\begin{equation}
\label{eq:33}
s^Tr\leq(\textbf{1}^Ts)f^{-1}\left[ \frac{\sum s_{(i)} f(r^{+}_{(i)})}{\textbf{1}^Ts}\right] 
\end{equation}
For concave nondecreasing invertible functions (for example $f(x)=\log_y(x)$, $y\geq1$) a similar result is obtainable in which the lower bound of $r$ should be used, as $r\succeq r^{-}$. If only the direction of $s=\|s\|_2e_s$ ($\|e_s\|_2=1$) is important, without loss of generality, $s$ can be defined as $s=(\sum e_{_i})e_s$, so we can eliminate the terms of $\textbf{1}^Ts$ in the previous inequalities. Therefore, \eqref{eq:32} can be reduced to
\begin{equation}
\label{eq:34}
s^Tr\leq f^{-1}\left[\sum s_{(i)} f(r_{(i)})\right] 
\end{equation}

Considering $f(x)=\|x\|$ as a convex function, Lemma \ref{lem5} leads to the following theorem:

\begin{theorem}
\label{th3}
Let $s\in\mathbb{R}^3$ be a unit vector such that $e_s\succeq0$. Assume that in the time interval of $\mathcal{T}$, it is known that the 2-norm of the CS's trajectory is bounded below such that $\rho^{-}=\inf_{t\in\mathcal{T}}(\|r\|_2)$. Moreover, assume that the maximum distance of an element of $r$ from the origin is bounded above by $\rho^{+}_{(i)}=\sup_{t\in\mathcal{T}}(\|r_{(i)}\|_1)$. Then, $\theta=\angle(s,r)$ is restricted by the following inequality
\begin{equation}
\label{eq:35}
\|\cos\theta\|_1\leq\min\left\{1,\frac{e_s^T\rho^{+}}{\rho^{-}}\right\}
\end{equation}
\end{theorem}
\begin{IEEEproof}
From the result of Lemma \ref{lem5}, considering $f$ to be any norm function $f(x)=\norm{x}$, since $x$ is a scalar variable, the function equivalently reduces to $f(x)=\norm{x}_1$. Thus, we have
\begin{equation}
\label{eq:36}
\norm{\frac{e_s^Tr}{\textbf{1}^Te_s}}_1 \leq\frac{\sum {e_s}_{(i)} \|r_{(i)}\|_1}{\textbf{1}^Te_s}
\end{equation}
The denominators are non-negative ($\textbf{1}^Te_s\geq0$), so \eqref{eq:36} reduces to
\begin{equation}
\label{eq:37}
\|e_s^Tr\|_1 \leq\sum {e_s}_{(i)} \|r_{(i)}\|_1
\end{equation}
The left-hand side of \eqref{eq:37} is equal to $\|e_s^Tr\|_1=\|r\|_2\|\cos\theta\|_1$. Therefore, taking $\|r\|_2$ into the denominator of the right-hand side of \eqref{eq:37} and considering the upper bound of the right-hand side by replacing the supremum of nominator and the infimum of denominator over time, inequality \eqref{eq:35} is proved. 
\end{IEEEproof}

\textcolor{black}{Conic bounds introduced in Theorem \ref{th3} particularly may be applicable in spacecraft formation flying as restricts the CS in a special TS's field of view in which can be considered as a beneficial property for missions equipped by vision-based relative navigation sensors.}

\begin{remark}
\label{rem7}
Theorem \ref{th3} defines a restricting cone which is also a function of $e_s$. In order to obtain the smallest set for $\theta$, the vector $e_s$ in the inequality \eqref{eq:35} should be selected in order to minimize the term $\sum {e_s}_{(i)}\rho^{+}_{(i)}$. Let $i^*$ be the nontrivial solution of $\min_i(\rho^{+}_{(i)})=\min_i[\sup_{t\in\mathcal{T}}(\|r_{(i)}\|_1)]$ (supposing $\rho^{+}_{(i^*)}\neq0$). Then, the vector $e_s$ that is constructed as ${e_s}_{(i=i^*)}=1$ and ${e_s}_{(i\neq i^*)}=0$ yields the smallest set for $\theta$ (and consequently for $r$).
\end{remark}

\subsection{Impulse Design Under Path Constraints}
\label{subsec:III-C}

In this subsection it is shown that how the impulse positions and times can be determined in order to satisfy the PCs subject to the CW system. Some boundaries are introduced regarding the Jordan-Brouwer seperation theorem, in which the satisfaction of the constraints is guaranteed by considering those boundary sets. 

The two-impulse trajectory of the CS subject to the initial and final positions is determined by two parameters of $t$ and $\Delta t_{i+1,i}$ that is formulated in \eqref{eq:9}. The set $\mathcal{R}_{i,i+1}(t)$ contains all position vectors which can be reached at $t\in\mathcal{T}_{\Delta t_{i+1,i}}$ by any $\Delta t_{i+1,i}\in\mathcal{T}_{\pi/\kappa}$. Now, consider the following set which is used further: 
\begin{equation}
\label{eq:38}
\mathcal{Q}_{i,i+1}=\bigcup_{t\in\mathcal{T}_{\pi/\kappa}}\mathcal{R}_{i,i+1}(t)
\end{equation}

Set $\mathcal{Q}_{i,i+1}$ is a three dimensional surface in which subject to $r=r_i$ at $t_i$ and $r=r_{i+1}$ at $t_{i+1}$ as the trivial constraints (using two impulses), $r\in\mathcal{Q}_{i,i+1}$ at $\forall t\in\Delta t_{i+1,i}$. Roughly speaking, the CS's trajectory entirely lies in the set $\mathcal{Q}_{i,i+1}$. 

\begin{lemma}[Jordan-Brouwer]
\label{lem6}
Let $\partial\mathcal{S}$ be a connected surface that is closed as a subset of $\mathbb{R}^3$. Then $\mathbb{R}^3-\partial\mathcal{S}$ has exactly two connected components ($\mathcal{S}$ and its complement $\mathbb{R}^3\setminus\mathcal{S}$) whose common boundary is $\partial\mathcal{S}$. 
\end{lemma}
\begin{IEEEproof}
See proof of Theorem 4.16 in \cite{b4}. 
\end{IEEEproof}

A result of Lemma \ref{lem6} is that every continuous path connecting a point in $\mathcal{S}$ to a point in $\mathbb{R}^3\setminus\mathcal{S}$ intersects somewhere with $\partial\mathcal{S}$. Thus, the following proposition can be directly concluded. 
\begin{proposition}
\label{prop2}
Let $t\in\mathcal{T}_{\pi/\kappa}$ and $r_s\in\mathbb{R}^3$ be reachable in $\Gamma_i^j(2)$ at a $t\in\mathcal{T}_{\pi/\kappa}$. Let $\mathcal{S}\subset\mathbb{R}^3$ be a set of position vectors with $\partial\mathcal{S}\subset\mathbb{R}^3$ as a boundary such that $\partial\mathcal{S}$ is a connected closed surface and $r_s\notin\mathcal{S}$. If the boundary of $\mathcal{S}$ is unreachable in $\Gamma_i^j(2)$ at any $t\in\mathcal{T}_{\pi/\kappa}$, then any member of $\mathcal{S}$ is unreachable in $\Gamma_i^j(2)$ at any $t\in\mathcal{T}_{\pi/\kappa}$, i.e.:   
\begin{equation}
\label{eq:39}
\partial\mathcal{S}\cap\mathcal{Q}_{i,i+1}=\varnothing\Rightarrow\mathcal{S}\cap\mathcal{Q}_{i,i+1}=\varnothing
\end{equation}
\end{proposition}
\begin{IEEEproof}
The set $\mathcal{Q}_{i,j}$ is a continuous surface in which the assumption $r_s\notin\mathcal{S}$ states that the surface has at least one point in $\mathbb{R}^3\setminus\mathcal{S}$. So, according to Lemma \ref{lem6}, if there exists a $t\in\mathcal{T}_{\pi/\kappa}$ in which $\partial\mathcal{S}\cap\mathcal{R}_{i,j}(t)\neq\varnothing$ then the curve intersects by $\partial\mathcal{S}$ somewhere. Therefore, its contrapositive is equivalently true. 
\end{IEEEproof}

The above-mentioned results can be subjected to the inequality PCs of type \eqref{eq:14} to conclude the following corollary. 

\begin{corollary}
\label{cor1}
Suppose Assumptions \ref{ass1} and \ref{ass2} hold. Assume it is known that a position vector $r_s\in\{r\in\mathbb{R}^3|\rho_l^{\prime}\leq\|r-\tilde{r}_l\|_2\leq \rho_l^{\prime\prime}\}$ is reachable in $\Gamma_i^j(2)$ at $\forall t\in\mathcal{T}_{\tilde{t}_l}$. An inequality PC of type \eqref{eq:14} is satisfied in $\Gamma_i^j(2)$ if any $r\in\{r\in\mathbb{R}^3|\|r-\tilde{r}_l\|_2=\rho_l^{\prime}\vee\|r-\tilde{r}_l\|_2=\rho_l^{\prime\prime}\}$ is unreachable in $\Gamma_i^j(2)$ at $\forall t\in\mathcal{T}_{\tilde{t}_l}$.
\end{corollary}
\begin{IEEEproof}
It can be simply proved by substituting $\mathcal{S}=\{r\in\mathbb{R}^3|\rho_l^{\prime}\leq\|r-\tilde{r}_l\|_2\leq \rho_l^{\prime\prime}\}$ and consequently $\partial\mathcal{S}=\{r\in\mathbb{R}^3|\|r-\tilde{r}_l\|_2=\rho_l^{\prime}\vee\|r-\tilde{r}_l\|_2=\rho_l^{\prime\prime}\}$ in Proposition \ref{prop2}. 
\end{IEEEproof}
\begin{remark}
\label{rem8}
A direct result from the above corollary can be presented for special kinds of inequality PCs. Suppose Assumptions \ref{ass1} and \ref{ass2} hold. Assume it is known that a position vector $r_s\in\{r\in\mathbb{R}^3|\rho_l^{\prime}\leq\|r-\tilde{r}_l\|_2\}$ is reachable in $\Gamma_i^j(2)$ at $\forall t\in\mathcal{T}_{\tilde{t}_l}$. Consider an inequality PC of type \eqref{eq:14} such that $\rho_l^{\prime\prime}\rightarrow\infty$ and $\tilde{t}=\pi/\kappa$. Then, the PC is satisfied in $\Gamma_i^j(2)$ if $\forall r\in\{r\in\mathbb{R}^3|\|r-\tilde{r}_l\|_2=\rho_l^{\prime}\}$, $r$ is unreachable in $\Gamma_i^j(2)$ at $\forall t\in\mathcal{T}_{\tilde{t}_l}$. 
\end{remark}
\begin{remark}
\label{rem9}
As an example, the vector $r_s$ that is used in Corollary \ref{cor1} can be selected to be the position at $t_1=0$, $r_s=r_1$, if $r_1$ is located inside $\widetilde{\mathcal{R}}_l$. This condition checks that if the CS is initially located inside or outside the set $\widetilde{\mathcal{R}}_l$. 
\end{remark}

Another result can be obtained from Lemma \ref{lem6} in which helps to find those impulse times (associated with fixed impulse positions) that the corresponding trajectory satisfies the constraints. 
\begin{proposition}
\label{prop3}
Suppose Assumptions \ref{ass1} and \ref{ass2} hold. Assume it is known that a point $r_s\in\mathcal{S}$ is reachable in $\Gamma_i^{i+1}(2)$ at $t_s\in[t_i,t_{i+1}]$ corresponding to $\Delta t_{i+1,i}=\Delta t_s$. Consider two time intervals of $\Delta t_a$ and $\Delta t_b$ such that $\Delta t_a<\Delta t_b$ and any point in $\mathcal{S}$ is unreachable in $\Gamma_i^{i+1}(2)$ at $\forall t\in[t_i,t_{i+1}]$ corresponding to both $\Delta t_{i+1,i}=\Delta t_a$ and $\Delta t_{i+1,i}=\Delta t_b$. Then, the following statements hold: 
\begin{enumerate}[label=\emph{(\roman*)}]
\item If $\Delta t_s\in[\Delta t_a,\Delta t_b]$ then any $r\in\mathcal{S}$ is unreachable in $\Gamma_i^{i+1}(2)$ at any $t\in[t_i,t_{i+1}]$ corresponding to $\Delta t_{i+1,i}\in[0,\Delta t_a]\cup[\Delta t_b,\pi/\kappa]$. 
\item If $\Delta t_s\in[0,\Delta t_a]\cup[\Delta t_b,\pi/\kappa]$ then any $r\in\mathcal{S}$ is unreachable in $\Gamma_i^{i+1}(2)$ at any $t\in[t_i,t_{i+1}]$ corresponding to $\Delta t_{i+1,i}\in[\Delta t_a,\Delta t_b]$. 
\end{enumerate}
\end{proposition}
\begin{IEEEproof}
Denote the set $\mathcal{R}_{i,i+1}(t)$ corresponding to $\Delta t_{i+1,i}=\Delta t_{s/a/b}$ by $\mathcal{R}_{s/a/b}(t)$. From Theorem \ref{th1} (the time uniqueness property), we know that $\mathcal{R}_{a}(t)$ and $\mathcal{R}_{b}(t)$ do not intersect, unless at $\Delta t_{i+1,i}=0$. The union of $\mathcal{R}_{a}(t)$ and $\mathcal{R}_{b}(t)$ can construct a boundary for a set of position vectors in $\mathcal{Q}_{i,i+1}$ which we refer to it by $\partial\mathcal{\widetilde{R}}_{a\cup b}(t)\subset\mathcal{Q}_{i,i+1}$ (and its interior by $\mathcal{\widetilde{R}}_{a\cup b}\subset\mathcal{Q}_{i,i+1}$). According to Proposition \ref{prop2}, since $\mathcal{R}_{a}(t)$ and $\mathcal{R}_{b}(t)$ do not intersect by $\mathcal{S}$, the set $\mathcal{\widetilde{R}}_{a\cup b}$ have no intersections by $\mathcal{S}$, i.e., $\mathcal{\widetilde{R}}_{a\cup b}\cap\mathcal{S}=\varnothing$, if $r_s\notin\mathcal{\widetilde{R}}_{a\cup b}$ (and $\mathcal{S}\subset\mathcal{\widetilde{R}}_{a\cup b}$ if $r_s\in\mathcal{\widetilde{R}}_{a\cup b}$). Item (i) can be proved by considering the fact that if $\Delta t_a\leq\Delta t_s\leq\Delta t_b$, then $\mathcal{R}_{s}(t)\subset\mathcal{\widetilde{R}}_{a\cup b}$. Roughly speaking $\mathcal{\widetilde{R}}_{a\cup b}$ is the set of those trajectories that $\Delta t_{i+1,i}\in[\Delta t_a,\Delta t_b]$. Item (ii) can be proved similarly by considering the fact that if $\Delta t_s\leq\Delta t_a$ or $\Delta t_b\leq\Delta t_s$, then $\mathcal{R}_{s}(t)\cap\mathcal{\widetilde{R}}_{a\cup b}=\varnothing$.
\end{IEEEproof}

It is worth mentioning that an $n$-impulse mission can be divided into a number of $n-1$ two-impulse missions and consequently any result about the two-impulse trajectories can be used for $n$-impulse cases just by reusing the result for $n-1$ times.

\section{Numerical Examples}
\label{sec:IV}

\subsection{Approximate Circular Formation Keeping}
\label{subsec:IV-A}

In this subsection the impulsive approximate CFK problem is analyzed by the use of numerical analysis. The CFK problem seeks for control solutions in order to keep the CS on a circular path around the TS. The term ``approximate'' is used to show that the CS is not going to lie on an exact circle (since is impossible with finite number of impulses), but should lie in a ring that is restricted by two circles; the approach circle ($\|r\|_2=\rho^{\prime\prime}$) and the keep-out circle ($\|r\|_2=\rho^{\prime}$).

\begin{figure}[!h]
\centering\includegraphics[width=0.5\linewidth]{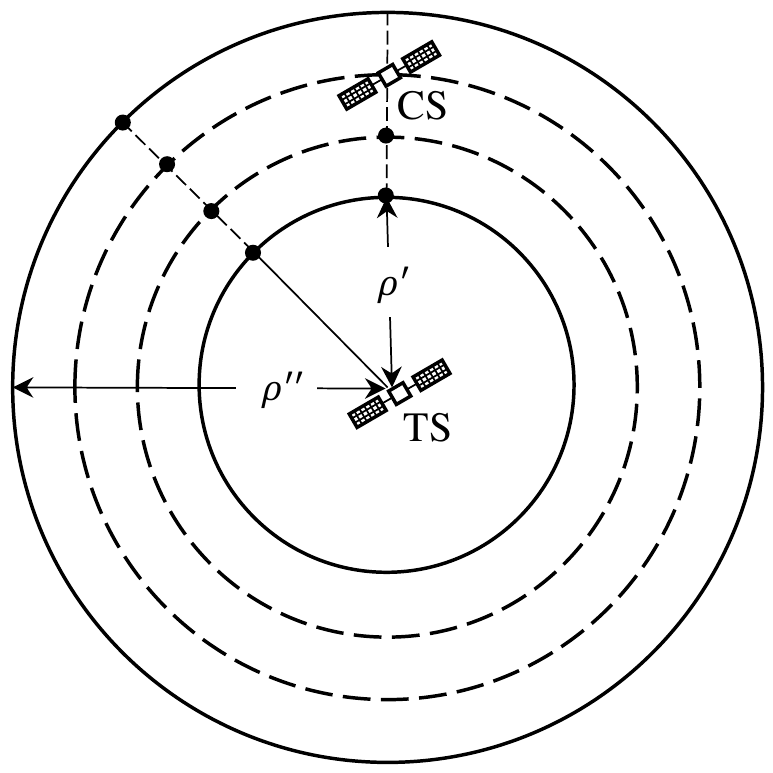}
\caption{Schematic grids for numerical analysis of impulsive approximate CFK.}
\label{fig:4}
\end{figure}

A finite number of polar grids are used for numerical computations, see Fig. \ref{fig:4}. The time is also approximated by discrete time instances. The numerical analysis decreases the accuracy of results based on how many nodes are used. However, increasing the number of nodes confronts the problem with the curse of dimensionality.

\begin{figure}[!h]
\centering\includegraphics[width=1\linewidth]{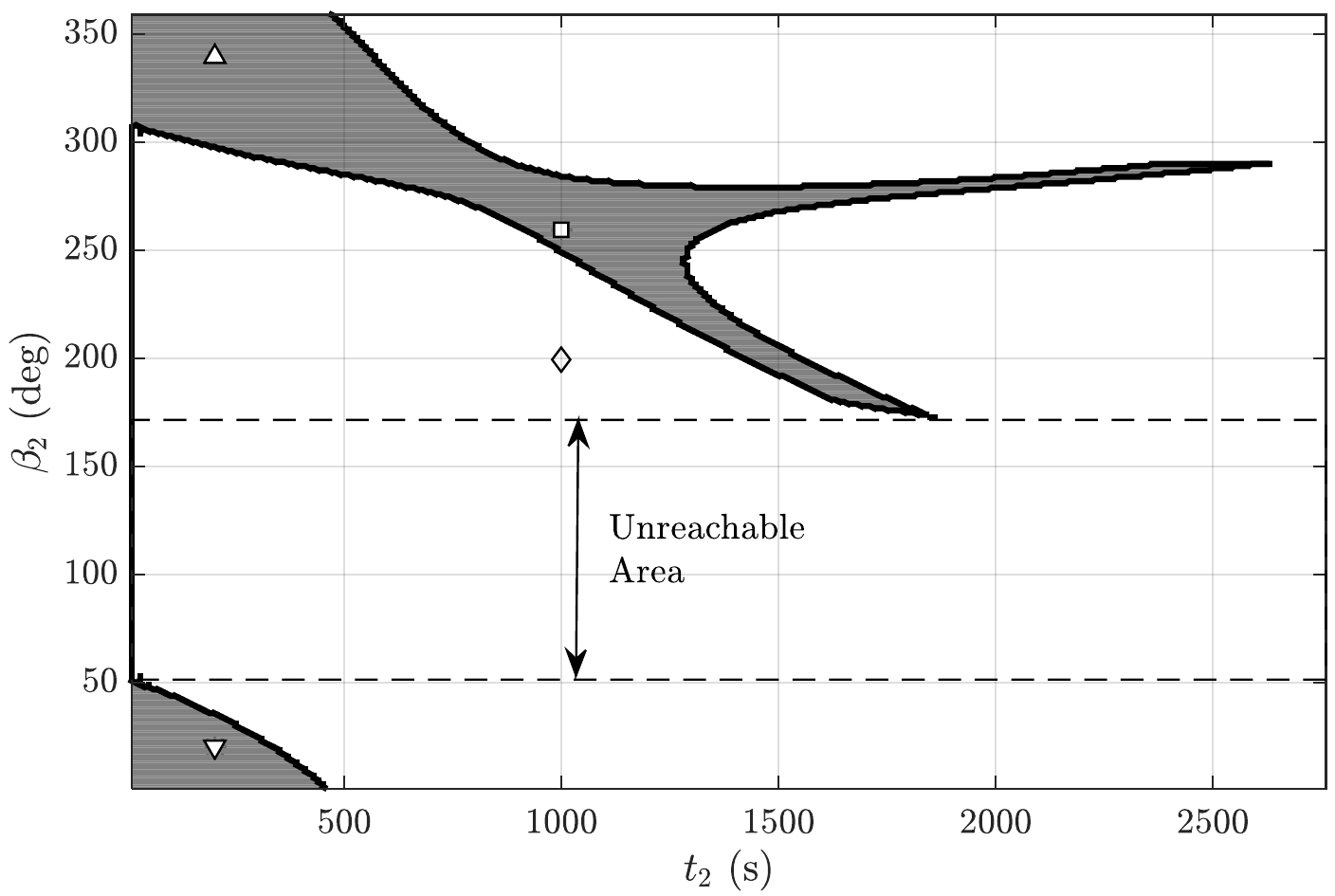}
\caption{Reachable region for the second impulse (highlighted in gray) and the unreachable area.}
\label{fig:5}
\end{figure}

\begin{figure}[!h]
\centering\includegraphics[width=1\linewidth]{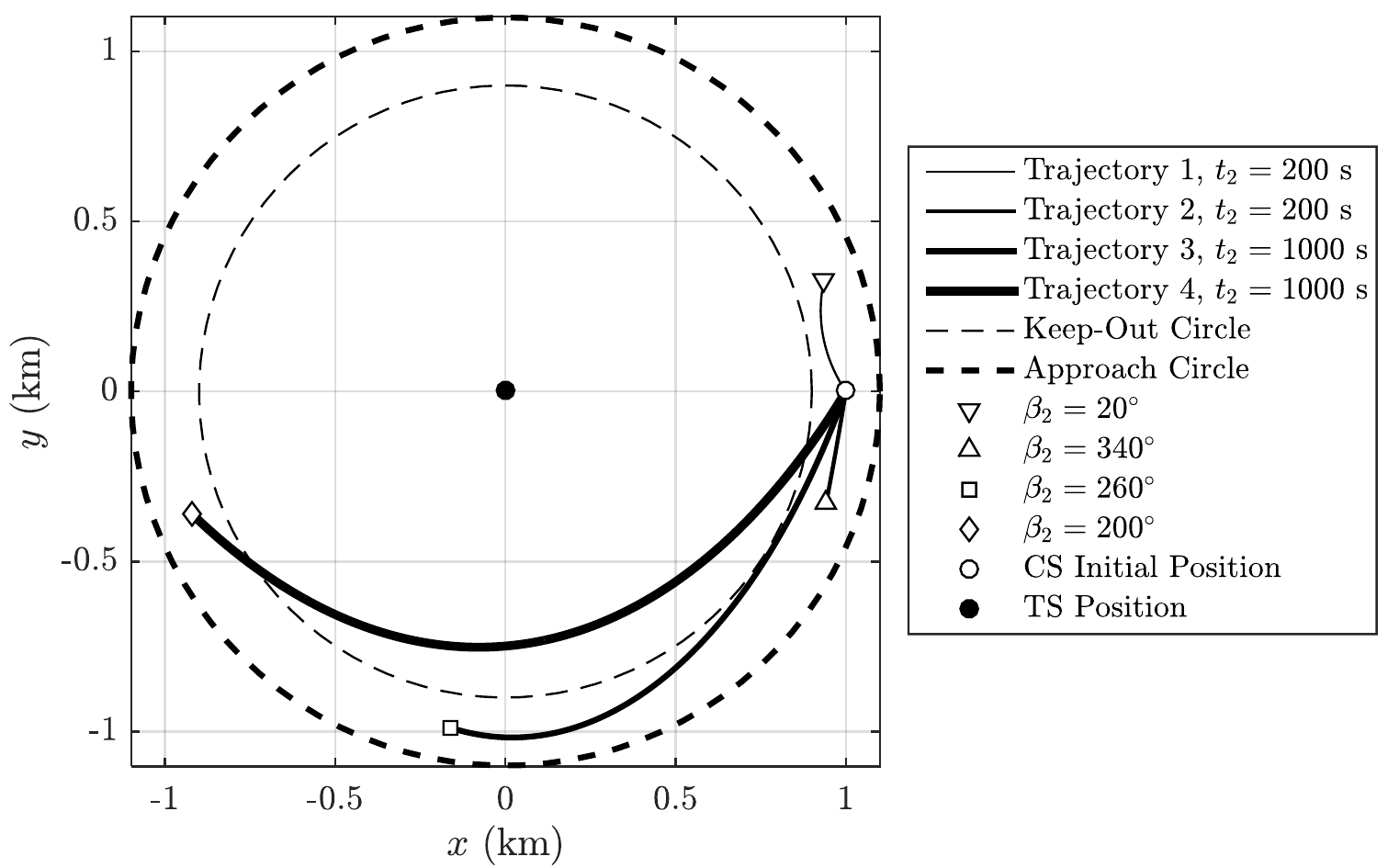}
\caption{Geometry of four sample trajectories from the first to the second impulse and the problem constraints.}
\label{fig:6}
\end{figure}

A two dimensional problem in $x$-$y$ plane is used as an example. The time step in the simulation is $10 \text{ s}$ and the TS is located at a circular orbit with an altitude of $400\text{ km}$. The radius of the approach and the keep-out circles are $\rho^{\prime}=0.9\text{ km}$ and $\rho^{\prime\prime}=1.1\text{ km}$, respectively. The position vectors for the impulses are considered to be $r_i=[\cos\beta_i\quad \sin\beta_i\quad 0]^T$, $i=2,3,...,n-1$. It is assumed that $\beta_i\in\mathbb{N}$, $0^\circ\leq\beta_i\leq360^\circ$, and $\beta_1=\beta_n=0$. In Fig. \ref{fig:5} those values of $\beta_2$ and $t_2$ in which the CS's trajectory in $\Gamma_1^2(2)$ satisfies the constrained problem is highlighted in gray. The area between dashed lines contains those positions that are unreachable in $\Gamma_1^2(2)$ at any $t\in\mathcal{T}_{t_2}$ corresponding to any $t_2\in\mathcal{T}_{\pi/\kappa}$ subject to the problem constraints. Therefore, the positions that are located between the dashed lines cannot be reached by two impulses and may become reachable by three impulses or more. Fig. \ref{fig:6} shows Four typical trajectories corresponding to the points that are specified in Fig. \ref{fig:5}. 

\begin{figure}[!h]
\centering\includegraphics[width=1\linewidth]{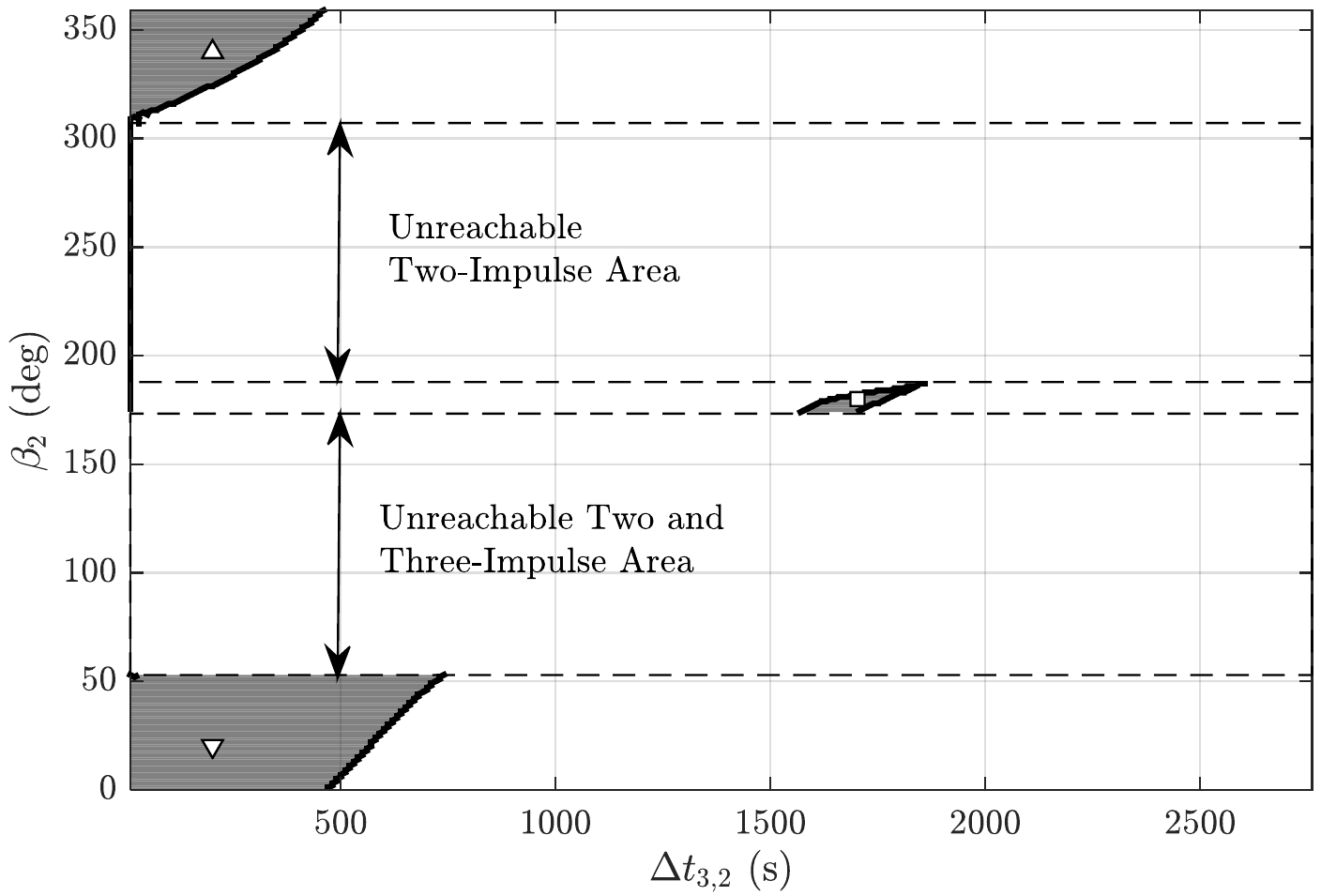}
\caption{Valid region for the second impulse (highlighted in gray) and the unreachable area.}
\label{fig:7}
\end{figure}

Fig. \ref{fig:7} shows those values of $\beta_2$ that are reachable from $i=1$ such that the $\beta_3=\beta_n=0$ is reachable from $i=2$ subject to the problem constraints. The highlighted areas of Fig. \ref{fig:7} can be divided into two categories. Consider the first category to be the union of two separated areas with $\Delta t_{3,2}\leq740\text{ s}$ at the left side (up and down) of the figure, and the second category to be the area with $1570\text{ s}\leq\Delta t_{3,2}\leq1860\text{ s}$ at the middle of the figure. The first category contains those trajectories that satisfy the constraints of the problem while the CS do not visit all values of the polar angles $0^\circ\leq\beta\leq360^\circ$ with respect to the TS. The second category includes those trajectories that the CS visits all values of $\beta$ and satisfy the constraints of the problem. 

\begin{figure}[!h]
\centering\includegraphics[width=1\linewidth]{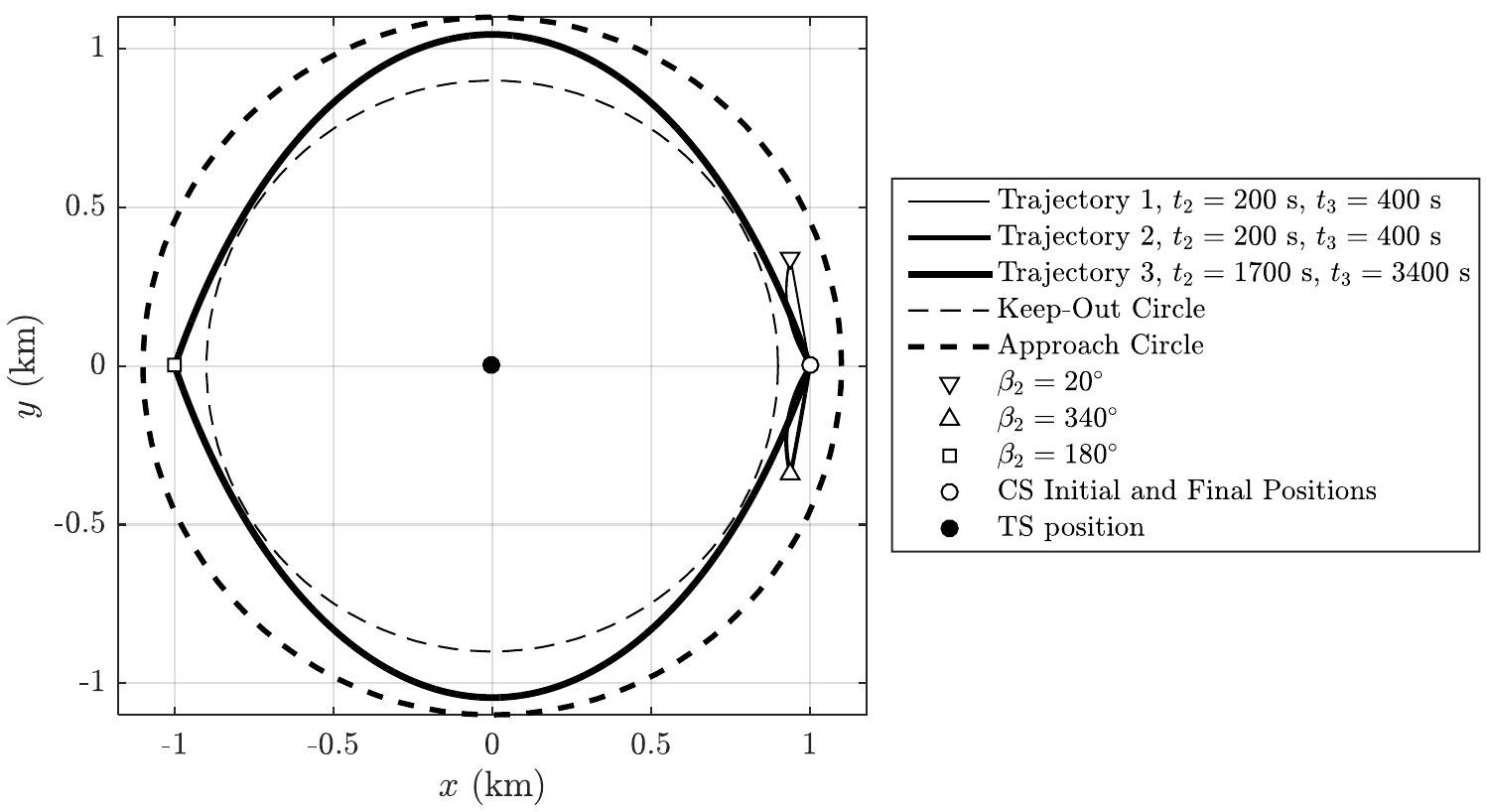}
\caption{Geometry of three sample trajectories from the first to the second and from the second to the third impulse in which the problem constraints are satisfied.}
\label{fig:8}
\end{figure}

Two unreachable areas are distinguished in Fig. \ref{fig:7}; the unreachable two- and three-impulse area, and the unreachable two- impulse area. The former includes the same unreachable points that are defied in Fig. \ref{fig:7}, and the latter contains the new results. The two and three-impulse unreachable points are those values of $\beta_2$ that cannot be reached from $i=1$. The two-impulse unreachable area are those values of $\beta_2$ that $\beta_3=\beta_n=0$ becomes unreachable from $i=2$. Fig. \ref{fig:8} shows three typical trajectories corresponding to the points that are specified in Fig. \ref{fig:7}. 

\subsection{Collision-Free Maneuver}
\label{subsec:IV-B}

In this subsection the impulsive CFM is analyzed numerically. The CFM is achieved by implementing those impulse positions at $i$ and $i+1$ such that a sphere area (problem constraint) do not be violated for every value of $\Delta t_{i+1,i}$. In this example fixed impulse positions is found such that the CS accomplish its mission while the constraints be satisfied independent from the transfer times. From Propositions \ref{prop2} and \ref{prop3}, we know that if $r_i$ and $r_{i+1}$ be considered such that the CS's trajectory corresponding to $\Delta t_{i+1,i}=0$ and $\Delta t_{i+1,i}=\pi/\kappa$ do not collide with the the constraint sphere, then the CFM is solved as well.

The trajectory of the CS corresponding to $\Delta t_{i+1,i}=0$ is a straight line from $r_i$ to $r_{i+1}$. The trajectory of the CS assuming $\Delta t_{i+1,i}=\pi/\kappa$ becomes singular, therefore, we approximate its locus by $\Delta t_{i+1,i}=\pi/\kappa-\epsilon$, such that $\epsilon>0$ is a small value to be determined. 

\begin{figure}[!h]
\centering\includegraphics[width=1\linewidth]{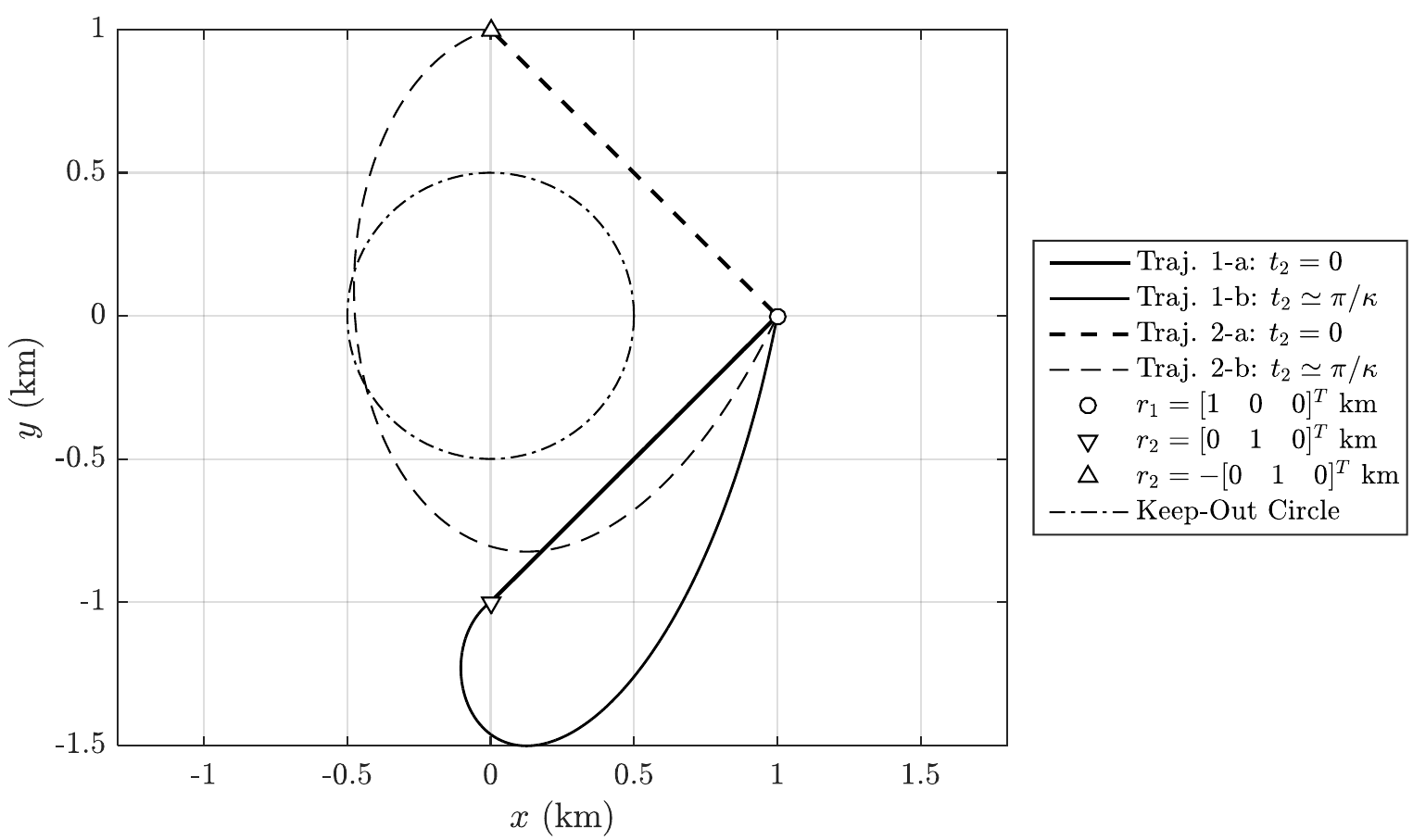}
\caption{Two choices for $r_2$ and the corresponding trajectory sets which are bounded by Traj. $i$-a and Traj. $i$-b.}
\label{fig:9}
\end{figure}

As an example consider a two dimensional problem in $x$-$y$ plane. Suppose the CS is located initially at $r_1=[1\quad0\quad0]^T \text{ km}$. The problem asks to find those values of $r_i$, $i\geq2$, in which a CFM can be accomplished such that the CS observes all the polar angles with respect to the TS and a keep-out circle with a radius of $\rho^{\prime}=0.5\text{ km}$ be satisfied. Assuming $\epsilon=1\text{ s}$, Fig. \ref{fig:9} shows two different chooses of $r_2=-[0\quad1\quad0]^T\text{ km}$ and $r_2=[0\quad1\quad0]^T\text{ km}$ beside the reachable area in $\Gamma_1^2(2)$ starting from $r_1$ and ending in $r_2$. In Fig. \ref{fig:9} it is shown that $r_2=-[0\quad1\quad0]^T\text{ km}$ leads to a CFM since both trajectories corresponding to $\Delta t_{i+1,i}=0$ and $\Delta t_{i+1,i}=\pi/\kappa-\epsilon\simeq\pi/\kappa$ do not collide with the keep-out circle. Instead, selecting the second impulse position to be $r_2=[0\quad1\quad0]^T\text{ km}$ the CFM cannot be constructed. 

\begin{figure}[!h]
\centering\includegraphics[width=1\linewidth]{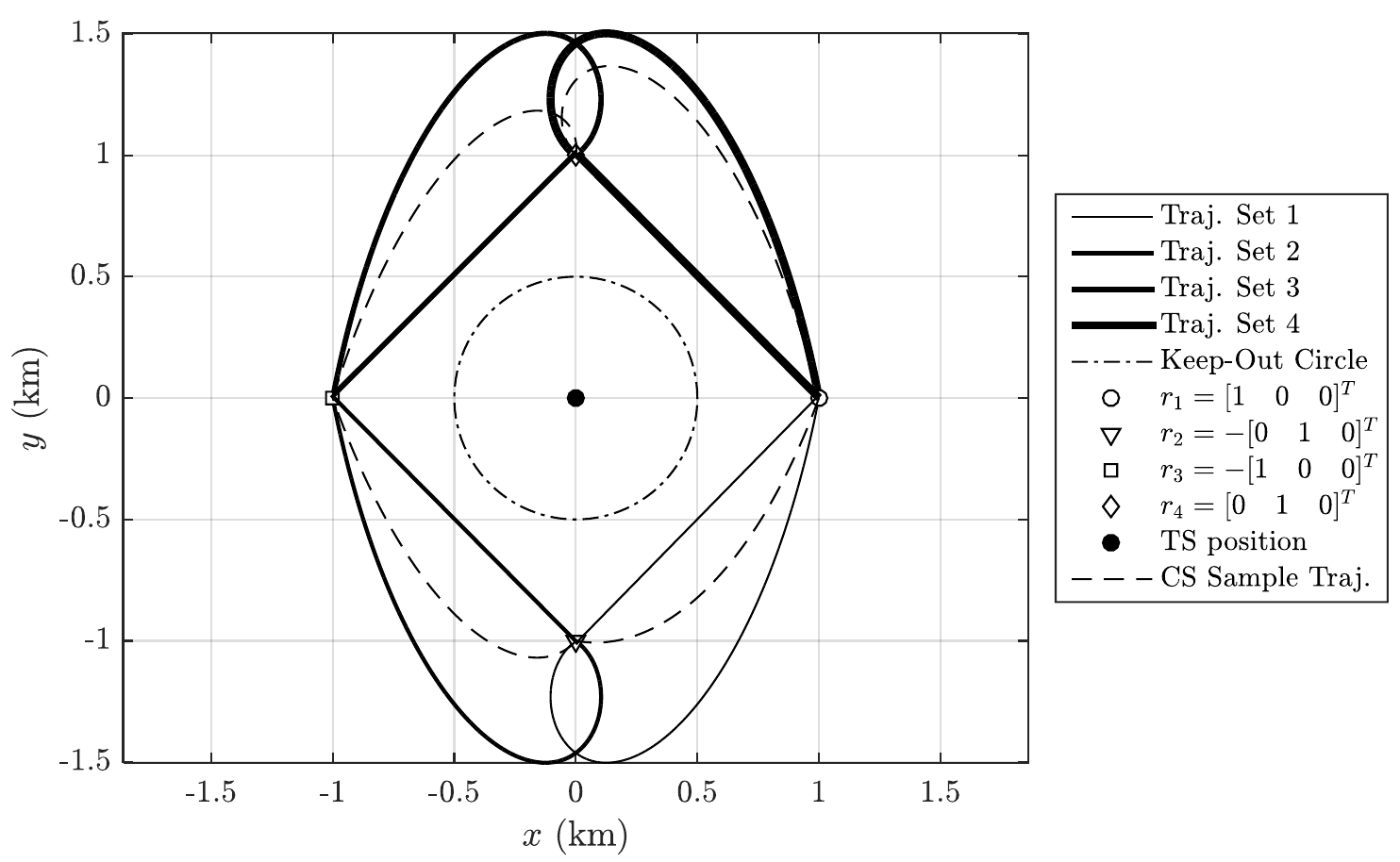}
\caption{Four impulse positions which are a solution to the CFM problem.}
\label{fig:10}
\end{figure}

Fig. \ref{fig:10} shows four impulse positions that can be considered as a solution to our CFM problem. Taking these four positions, for any time intervals $\Delta t_{i+1,i}\in\mathcal{T}_{\pi/\kappa}$, $i=1,2,3$, the keep-out circle is not violated and a periodic motion around the TS can be accomplished as well. 

\subsection{Discussions}
\label{sec:subsec:IV-C}

In this section some relations between Theorems \ref{th2} and \ref{th3} with the numerical examples of Section \ref{sec:IV} are discussed in more detail. 

Theorem \ref{th2} introduces an upper norm bound for CS's trajectory. This upper norm bound is tabulated for every trajectory example of Section \ref{sec:IV} in Table \ref{table:1}. In Figs. \ref{fig:6}, \ref{fig:8}, \ref{fig:9}, and \ref{fig:10}, the impulse positions are located on a unit circle, therefore if $\Delta t_{i+1,i}<(\pi/\kappa)/2\simeq1428\text{ s}$ then the upper bound (according to Theorem \ref{th2}) is simply $\sqrt{2}$ (such as N$\underbar{o}$s. 1-8 and 11-14 in Table \ref{table:1}). If $\Delta t_{i+1,i}>(\pi/\kappa)/2\simeq1428\text{ s}$ then the upper bound exceeds $\sqrt{2}$ and is computed by \eqref{eq:23} and \eqref{eq:24} (such as N$\underbar{o}$s. 9 and 10 in Table \ref{table:1}). 

Theorem \ref{th3} introduces conic bounds on the CS's trajectory. Using the optimal index value $i^*$ (which is discussed in Remark \ref{rem7}), the conic bounds for every trajectory example of Section \ref{sec:IV} is presented in Table \ref{table:2}. Each no. corresponds to a two-impulse trajectory which is previously defined in Table \ref{table:1}. For N$\underbar{o}$s. 3, 4, and 9-14, the upper bound of $\|\cos\theta\|_1$ (i.e., $c_{\theta}$) are equal to 1 that are obtained using \eqref{eq:35}. Therefore, the conic bound of Theorem \ref{th3} is useless for these scenarios. For N$\underbar{o}$s. 1, 2, and 5-8, the upper bound has a value less than unity which restricts the CS's trajectory to lie outside a double cone with an obtuse aperture. 

\begin{table}[ht]
\caption{The upper norm bounds, $\delta_{i+1,i}$ (obtained from Theorem \ref{th2}), for the two-impulse numerical examples of Section \ref{sec:IV}.}
\begin{center}
\begin{tabular}{|c|c|c|c|c|c|}
	\hline
	\multirow{2}{*}{N$\underbar{o}$.}&\multicolumn{4}{c|}{Example} &\multirow{2}{*}{$\delta_{i+1,i}$ (km)}\\
	\cline{2-5}
	\multicolumn{1}{|c|}{}&Fig.&$\beta_i$ (deg.)&$\beta_{i+1}$ (deg.)&$\Delta t_{i+1,i}$ (s)&\multicolumn{1}{c|}{}\\
	\hline
	1&\multirow{4}{*}{\ref{fig:6}}&\multirow{4}{*}{$0$}&$20$&\multirow{2}{*}{$200$}&\multirow{8}{*}{$\sqrt{2}$}\\ 
	\cline{1-1}\cline{4-4}
	2&\multicolumn{1}{c|}{}&\multicolumn{1}{c|}{}&$340$&\multicolumn{1}{c|}{}&\multicolumn{1}{c|}{}\\ 
	\cline{1-1}\cline{4-5}
	3&\multicolumn{1}{c|}{}&\multicolumn{1}{c|}{}&$260$&\multirow{2}{*}{$1000$}&\multicolumn{1}{c|}{}\\ 
	\cline{1-1}\cline{4-4}
	4&\multicolumn{1}{c|}{}&\multicolumn{1}{c|}{}&$200$&\multicolumn{1}{c|}{}&\multicolumn{1}{c|}{}\\
	\cline{1-1}\cline{2-5}
	5&\multirow{6}{*}{\ref{fig:8}}&\multirow{2}{*}{$0$}&$20$&\multirow{4}{*}{$200$}&\multicolumn{1}{c|}{}\\ 
	\cline{1-1}\cline{4-4}
	6&\multicolumn{1}{c|}{}&\multicolumn{1}{c|}{}&$340$&\multicolumn{1}{c|}{}&\multicolumn{1}{c|}{}\\ 
	\cline{1-1}\cline{3-4}
	7&\multicolumn{1}{c|}{}&$20$&\multirow{2}{*}{$0$}&\multicolumn{1}{c|}{}&\multicolumn{1}{c|}{}\\ 
	\cline{1-1}\cline{3-3}
	8&\multicolumn{1}{c|}{}&$340$&\multicolumn{1}{c|}{}&\multicolumn{1}{c|}{}&\multicolumn{1}{c|}{}\\ 
	\cline{1-1}\cline{3-6}
	9&\multicolumn{1}{c|}{}&$0$&$180$&\multirow{2}{*}{$1700$}&\multirow{2}{*}{$\simeq1.19\sqrt{2}$}\\ 
	\cline{1-1}\cline{3-4}
	10&\multicolumn{1}{c|}{}&$180$&$0$&\multicolumn{1}{c|}{}&\multicolumn{1}{c|}{}\\ 
	\hline
	11&\multirow{4}{*}{\ref{fig:9},\ref{fig:10}}&$0$&$270$&\multirow{2}{*}{$<(\pi/\kappa)/2$}&\multirow{4}{*}{$\sqrt{2}$}\\ 
	\cline{1-1}\cline{3-4}
	12&\multicolumn{1}{c|}{}&$270$&$180$&\multicolumn{1}{c|}{}&\multicolumn{1}{c|}{}\\ 
	\cline{1-1}\cline{3-4}
	13&\multicolumn{1}{c|}{}&$180$&$90$&\multirow{2}{*}{$\lesssim1428$}&\multicolumn{1}{c|}{}\\ 
	\cline{1-1}\cline{3-4}
	14&\multicolumn{1}{c|}{}&$90$&$0$&\multicolumn{1}{c|}{}&\multicolumn{1}{c|}{}\\ 
	\hline
	\multicolumn{6}{l}{\scriptsize{\textbf{Note:} $r_i=[\cos\beta_i\quad\sin\beta_i\quad0]^T$}}
\end{tabular}
\end{center}
\label{table:1}
\end{table}

\begin{table}[ht]
\caption{Conic bounds, $\|\cos\theta\|_1\leq c_{\theta}$ (obtained from Theorem \ref{th3}), for the two-impulse numerical examples of Section \ref{sec:IV}.}
\begin{center}
\begin{tabular}{|c|c|c|c|c|}
	\hline
	\multirow{2}{*}{N$\underbar{o}$.}&\multicolumn{3}{c|}{Example} & \multirow{2}{*}{$c_{\theta}$}\\
	\cline{2-4}
	\multicolumn{1}{|c|}{} & $e_s^T$ (deg.)&$\widehat{\rho}^-$ (km)&$(\widehat{\rho}^+)^T$ (km)&\multicolumn{1}{c|}{}\\
	\hline
	1 & \multirow{3}{*}{$[0\quad1\quad0]$} & \multirow{2}{*}{$0.9$} & \multirow{2}{*}{$[1\quad0.5\quad0]$} & \multirow{2}{*}{$5/9$} \\ 
	\cline{1-1}
	2 & \multicolumn{1}{c|}{} & \multicolumn{1}{c|}{} & \multicolumn{1}{c|}{} & \multicolumn{1}{c|}{} \\ 
	\cline{1-1}\cline{3-5}
	3 & \multicolumn{1}{c|}{} & $0.5$ & $[1\quad0.9\quad0]$ & \multirow{2}{*}{$1$} \\ 
	\cline{1-4}
	4 & $[1\quad0\quad0]$ & \multirow{7}{*}{$0.9$} & $[1\quad1.1\quad0]$ & \multicolumn{1}{c|}{} \\ 
	\cline{1-2}\cline{4-5}
	5 & \multirow{4}{*}{$[0\quad1\quad0]$} & \multicolumn{1}{c|}{} & \multirow{4}{*}{$[1\quad0.5\quad0]$} & \multirow{4}{*}{$5/9$} \\ 
	\cline{1-1}
	6 & \multicolumn{1}{c|}{} & \multicolumn{1}{c|}{} & \multicolumn{1}{c|}{} & \multicolumn{1}{c|}{} \\ 
	\cline{1-1}
	7 & \multicolumn{1}{c|}{} & \multicolumn{1}{c|}{} & \multicolumn{1}{c|}{} & \multicolumn{1}{c|}{} \\ 
	\cline{1-1}
	8 & \multicolumn{1}{c|}{} & \multicolumn{1}{c|}{} & \multicolumn{1}{c|}{} & \multicolumn{1}{c|}{} \\ 
	\cline{1-2}\cline{4-5}
	9 & \multirow{6}{*}{$[1\quad0\quad0]$} & \multicolumn{1}{c|}{} & \multirow{2}{*}{$[1\quad1.1\quad0]$} & \multirow{6}{*}{$1$} \\ 
	\cline{1-1}
	10 & \multicolumn{1}{c|}{} & \multicolumn{1}{c|}{} & \multicolumn{1}{c|}{} & \multicolumn{1}{c|}{} \\ 
	\cline{1-1}\cline{3-4}
	11 & \multicolumn{1}{c|}{} & \multirow{4}{*}{$\leq1$} & \multirow{4}{*}{$[1\quad\geq1\quad0]$} & \multicolumn{1}{c|}{} \\ 
	\cline{1-1}
	12 & \multicolumn{1}{c|}{} & \multicolumn{1}{c|}{} & \multicolumn{1}{c|}{} & \multicolumn{1}{c|}{} \\ 
	\cline{1-1}
	13 & \multicolumn{1}{c|}{} & \multicolumn{1}{c|}{} & \multicolumn{1}{c|}{} & \multicolumn{1}{c|}{} \\ 
	\cline{1-1}
	14 & \multicolumn{1}{c|}{} & \multicolumn{1}{c|}{} & \multicolumn{1}{c|}{} & \multicolumn{1}{c|}{} \\ 
	\hline
	\multicolumn{5}{l}{\scriptsize{\textbf{Note:} $\widehat{\rho}^+$ and $\widehat{\rho}^-$ are the estimated values of $\rho^+$ and $\rho^-$.}}
\end{tabular}
\end{center}
\label{table:2}
\end{table}

\begin{figure}[!h]
\centering\includegraphics[width=1\linewidth]{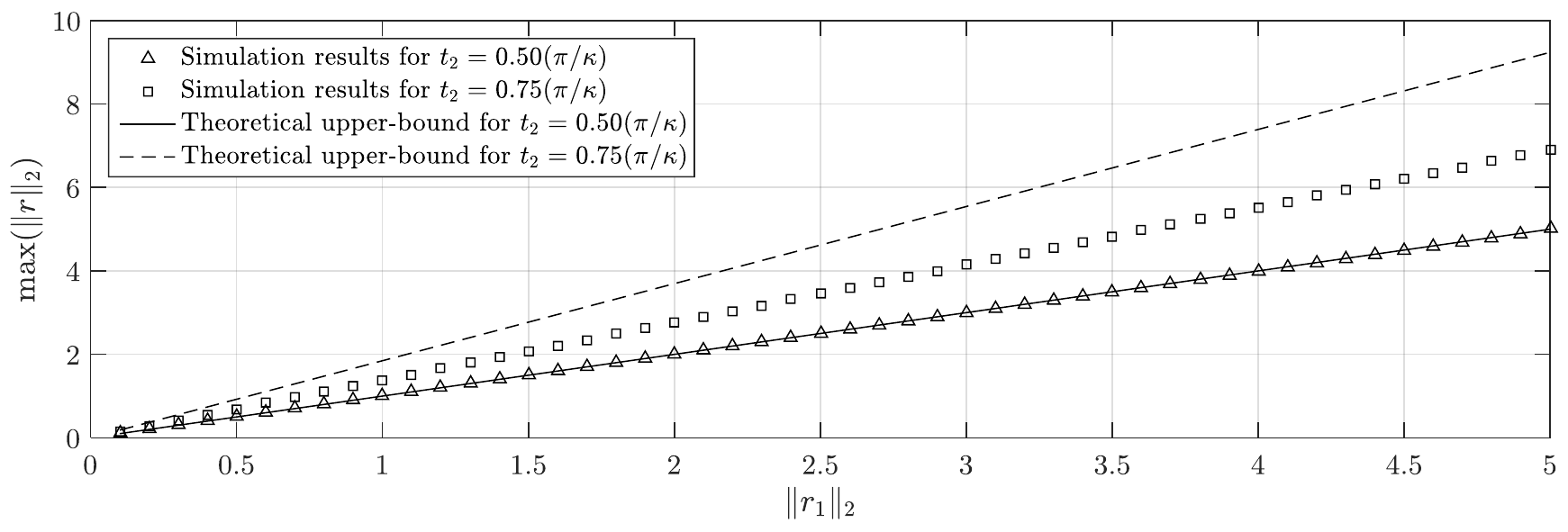}
\caption{Comparison of the maximum reached distance of the CS in the numerical simulation with the theoretical results as a function of the initial distance.}
\label{fig:11}
\end{figure}

In Fig. \ref{fig:11} a three dimensional example is used for the numerical evaluation of the trajectory upper-bound found in Theorem \ref{th2}. The TS is located at the origin with an altitude of $400\text{ km}$ above Earth. The initial location of the CS is formulated as $r_1=\|r_1\|_2[\cos(\phi)\cos(\psi)\quad\cos(\phi)\sin(\psi)\quad\sin(\phi)]^T$ in which the simulations are done for $\phi\in[-\pi,\pi]$, $\psi=[0,2\pi]$, and $t=[0,t_2]$, and the maximum reached distance defined as $\max_{\phi,\psi,t}(\|r\|_2)$ is evaluated for each amount of $\|r_1\|_2$ from $0.1$ to $5$. This method is accomplished separately for $t_2=0.5\pi/\kappa$ and $t_2=0.75\pi/\kappa$.

\section{Conclusions}
\label{sec:VI}

The relative spacecraft motion has been analyzed under path constraints using the Clohessy-Wiltshire (CW) equations. Initially, the time uniqueness of the spacecraft's trajectory is analyzed under assumptions and the main result is proved in a theorem. The spectral analysis of the CW equations demonstrates some facts which are used to determine upper norm bounds for the spacecraft position between adjacent impulses. Moreover, a finite form of the Jensen's inequality is implemented to develop a conic bound for the spacecraft path in which needs additional priory estimations about the position time-history. Furthermore, it is shown that the unreachability of a set of continuous position vectors can be proven by considering the unreachability of some boundary positions. Finally, two numerical examples in the $x$-$y$ plane are presented. The first example is an approximate circular formation keeping where seeks for those impulse positions and times such that the chaser spacecraft's trajectory lie in a ring. The second example is a collision-free maneuver in which the impulse positions are found such that the chaser spacecraft do not violate a keep-out circle with any choice of impulse times, i.e., attaining a set of trajectories that are robust in terms of impulse times.

\section*{Acknowledgment}

I would like to express my special thanks to Dr. Nima Assadian for his helpful advice.

\end{document}